\documentclass[12pt,a4paper]{article}
\usepackage[cp1251]{inputenc}
\usepackage{amsmath}
\usepackage{amssymb}
\usepackage{latexsym}
\usepackage{srcltx}

\usepackage{pdflscape}

\usepackage{subcaption}
\usepackage{tikz}
\usepackage{anyfontsize}
\usepackage{tabularx}

\usepackage{multirow}
\usepackage{array}

\usepackage{amssymb, amsmath, amsfonts, mathtools}
\usepackage{esint}

\usepackage{breqn}
\usepackage{float}
\usepackage{caption}
\usepackage{subcaption}
\usepackage{hhline}
\usepackage{multirow}
\usepackage{rotating}
\usepackage{fixltx2e}

\usepackage{hyperref}

\newcommand{\R}{{\mathbb R}}
\newcommand{\re}{{\mathbb R}}
\newcommand{\n}{{\mathbb N}}

\newcommand{\z}{{\mathbb Z}}

\newcommand{\cT}{{\mathcal{T}}}
\newcommand{\cH}{{\mathcal{H}}}

\newcommand{\cR}{{\mathcal{R}}}

\newcommand{\by}{{\boldsymbol y}}
\newcommand{\bx}{{\boldsymbol x}}

\newcommand{\ba}{{\boldsymbol a}}
\newcommand{\bb}{{\boldsymbol b}}
\newcommand{\bc}{{\boldsymbol c}}
\newcommand{\bd}{{\boldsymbol d}}

\newcommand{\bk}{{\boldsymbol k}}

\newcommand{\be}{{\boldsymbol e}}
\newcommand{\bm}{{\boldsymbol m}}
\newcommand{\bxi}{{\boldsymbol \xi}}
\newcommand{\bq}{{\boldsymbol q}}

\newcommand{\bv}{{\boldsymbol v}}

\newcommand{\bh}{{\boldsymbol h}}

\newcommand{\nul}{{\boldsymbol 0}}

\renewcommand{\Re}{\mathop{\mathrm{Re}}\nolimits}

\newtheorem{theorem}{Theorem}
\newtheorem{prop}{Proposition}
\newtheorem{lemma}{Lemma}
\newtheorem{cor}{Corollary}
\newtheorem{remark}{Remark}

\newtheorem{defi}{Definition}
\newtheorem{conj}{Conjecture}

\newtheorem{ser}{Series}

\newcommand{\overbar}[1]{\mkern 1.5mu\overline{\mkern-1.5mu#1\mkern-1.5mu}\mkern 1.5mu}

\date{}

\author{
Vladimir Yu. Protasov
\thanks{DISIM, University of L'Aquila (Italy), Moscow State University, Department of Mechanics and Mathematics, {e-mail: \tt\small
v-protassov@yandex.ru}} , 
Tatyana Zaitseva
\thanks{Moscow Center for Fundamental and Applied Mathematics, Moscow State University, Department of Mechanics and Mathematics, {e-mail: \tt\small
zaitsevatanja@gmail.com}} 
}

\title{Self-affine 2-attractors and tiles
\thanks{
The first  author is supported  by the Russian
Foundation for Basic Research, projects  no. 19-04-01227 and 20-01-00469
}}

\begin{document}
\maketitle

\begin{abstract}
We study  two-digit attractors (2-attractors) in~$\re^d$ which are self-affine compact sets defined by two 
contraction affine mappings with the same linear part. They are widely  studied in the literature under various names:  twindragons, two-digit tiles, 2-reptiles, etc., due to many applications 
in approximation theory, in the  construction of multivariate Haar systems and other wavelet bases, in the discrete geometry, and in the number theory. We obtain a complete classification of isotropic  
2-attractors~in~$\re^d$ and show that they are all homeomorphic but not diffeomorphic. 
In the general, non-isotropic, case it is proved that a 2-attractor is uniquely defined, up to 
an affine similarity, by the spectrum of the dilation matrix. We estimate the number of different 2-attractors in~$\re^d$ by analysing integer unitary expanding polynomials with the free coefficient $\pm 2$. The total number of such polynomials is estimated by the Mahler measure.  We present several infinite series of such polynomials. For some of the 2-attractors, their 
H\"older exponents are found. Some of our results are extended to attractors with 
an arbitrary number of digits.

\bigskip

\smallskip

\begin{flushright}

\end{flushright}

\end{abstract}
\bigskip

\begin{center}
\textbf{1. Introduction}
\end{center}
\bigskip

Self-similar attractor in~$\re^d$ is a  compact set defined by an integer  $d\times d$ 
matrix~$M$ and by a finite set of integer vectors 
(``digits'')~$D \, = \, \{\bd_0, \ldots , \bd_{m-1}\}\subset \z^d$, where $m = |{\rm det}\, M|$.
The matrix $M$ is assumed to be {\em expanding}, i.e., 
the moduli of all its eigenvalues are larger than one. The digits~$\bd_i \in D$ 
are taken from different quotient classes  $\z^d/M\z^d$, 
which means  $\bd_{i} - \bd_j \, \notin \, M\z^d$, whenever $i\ne j$.  
\begin{defi}\label{d.attractor}   
A self-affine attractor corresponding to an integer expanding matrix~$M$ 
and to a set of digits~$D$ is the following set 
\begin{equation}\label{eq.G}
G \ = \ G(M, D)\quad = \quad \left\{\, 0.\ba_1 \ba_2 \ldots \ = \ 
\sum_{k=1}^{\infty} M^{-k}\ba_k  \quad : \quad
\ba_k \in D \, \right\}.
\end{equation}
\end{defi}
A self-affine attractor is a compact set with a nonempty interior, 
its Lebesgue measure $|G|$ is always a natural number. The integer shifts 
 $\{G+k\}_{k \in \z^d}$ cover the entire space 
 $\re^d$ in  $|G|$ layers, i.e., every point 
 $x\in \re^d$ apart from a null set belongs to exactly $|G|$ shifts. In the case 
  $|G| = 1$ the attractor is called  {\em tile} and its shifts 
  form a  {\em tiling} of the space, i.e., its partition in one layer, see~\cite{GH, LW97}.

Formula~(\ref{eq.G}) means that the attractor 
plays a role of a unit segment in the space $\re^d$ in the 
$M$-adic system with digits from~$D$. However, even on the line~$\re$, 
an attractor can be different from the segment. For example, in the triadic system   
with digits  $0,1$, and $2$ the set $G$ is a segment, but if the digit  
 $2$ is replaced by  $5$, then  $G$ becomes a fractal-like  set 
with an infinite number of connected components~\cite{Woj, CP}. 
Similarly to the segment, any attractor is self-affine.   
For an arbitrary $\ba \in \z^d$,  we denote by~$M_{\ba}$ the affine  operator 
$M_{\ba}\, \bx \, = \, M\bx - \ba, \, \bx \in \re^d$. 
It is easily shown that 
$\ G \, = \, \bigcup\limits_{\ba \in D} M_{\ba}^{\, -1} \, G$,  
and since the measure of each of the sets  $M_{\ba}^{\, -1} \, G$ is equal to 
 $m^{-1}\, |G|$, and their sum is equal to~$|G|$, it follows that all of those sets 
 have intersections of zero measure.  Therefore,  $G$ 
is a disjunct, up to a null set, union of shifts of the set~$M^{-1}G$, 
i.e.,  $G$ is indeed self-affine. 
That is why the characteristic function~$\varphi = \chi_G$
of an arbitrary attractor is a  {\em refinable function}, which satisfies the following 
{\em refinement equation}
\begin{equation}\label{eq.ref}
\varphi(\bx) \  = \  \sum_{\bk \in D} \varphi (M\bx -\bk)
\end{equation} 
If the attractor is a tile, then the function~$\varphi$
possesses an orthonormal integer shifts. Consequently, it generates a multiresolution analysis (MRA)   of the space~$L_2(\re^d)$. 
That is why tiles are applied for construction of wavelets systems, in particular, 
of Haar bases in~$L_2(\re^d)$, see, for instance~\cite{BS, CHM, GM, KPS, LW95,  Woj, Zakh}. 

An attractor is called  {\em isotropic}, if it is generated by an 
{\em isotropic matrix}~$M$, which is a matrix with all eigenvalues of equal moduli but 
without nontrivial Jordan blocks.  An isotropic matrix is similar to 
an orthogonal matrix multiplied by a number. Isotropic attractors are most popular in applications. 

Two-digit attractors   (2-attractors) are those for which 
${\rm det}\, M \, = \pm 2$. In this case $m = |D| = 2$ and the
$M$-adic system  in $\re^d$ is similar to the dyadic system on the real line. 
In particular, each 2-attractor~$G$ disintegrates into two equal copies
$M^{-1}(G+\bd_0)$ and $M^{-1}(G+\bd_1)$ and the refinement equation~(\ref{eq.ref})
has only two terms. 
The Haar system has a unique generating function $\psi(\bx) = \varphi(M\bx-\bd_0) - 
\varphi(M\bx-\bd_1)$ unlike the general case, where it has  $m-1$ 
generating functions. For wavelets systems,  the situation is the same: 
 in the two-digit case the system of wavelets is generated by one 
 wavelet function and the corresponding MRA has a dyadic structure as for 
 wavelets in~$L_2(\re)$. That is why wavelets systems generated by 2-attractors 
 are natural and convenient in use. 

Two-digit attractors are studied in an extensive literature, a brief 
overview is given in the end of this section. 
It is known that for every $d$, there exist finitely many different, up to an affine similarity,  2-attractors in~$\re^d$. In~$\re^2$ there are precisely three 
2-attractors: a square, a dragon, and a bear (in the literature they are also known as rectangle, twindragon and tame twindragon).  They are all isotropic.  In $\re^3$ 
there exist seven 2-attractors, and there is only one isotropic among them, which is a cube. 
\smallskip 

\textbf{The fundamental results}. We derive a classification and geometric and topologic properties of 2-attractors. Some of our results, if possible, will be extended to arbitrary number of digits. In the isotropic case, the classification problem is solved completely (Section~5). That classification turns out to be rather simple, which is a bit surprising. Namely, in the odd dimensions~$d=2k+1$, all isotropic 2-attractors are parallelepipeds 
 (Theorem~\ref{th.isotr-odd}), while in every even dimension~$d=2k$, there exist precisely three, up to affine similarity, 2-attractors.  These are the parallelepiped, 
the direct product of  $k$ (two-dimensional) dragons, and the direct product of   $k$ (two-dimensional) bears  (Theorem~\ref{th.isotr-even}). The proofs  are constructive and the matrices~$M$ of these 2-attractors are explicitly found. 

The obtained classification allows us to establish many properties of isotropic 2-attractors. 
In particular, we show that all isotropic 2-attractors in~$\re^d$ are homeomorphic to 
the 
$d$-dimensional ball, but not  $C^1$-diffeomorphic to each other. Moreover, all these three types  of isotropic 2-attractors in an even-dimensional space are not metrically equivalent to each other. This means that they cannot be mapped to each other by bi-Lipschitz maps 
 (Theorem~\ref{th.isotr-nehomog} in Section~6). Properties of convex hulls of 2-attractors are analysed in Section~7.  
It turns out that among all  2-attractors, there are exactly two ones with polyhedral 
convex hull. This is a parallelepiped, whose convex hull has, of course,   $2^d$ vertices, 
and  a dragon, whose convex hull has  
$2^{3d/2}$ vertices. For all other   2-attractors, the convex hull is a zonoid with infinitely many extreme points. 

 In Section~8 we consider applications. We show that every isotropic 2-attractor is 
 a tile and hence it generates an MRA and a Haar system in~$L_2(\re^d)$
 provided the matrix~$M$ and the set of digits~$D$ are chosen in an appropriate way
(Theorem~\ref{th.tile-iso} gives this choice explicitly). 
Let us emphasise  that the Haar function generated by a direct  product of the bivariate dragons is 
different from the product of bivariate Haar functions! 
Analogous statements for general (non-isotropic) 2-attractors are formulated in a conjecture,
which is verified in small dimensions (Proposition~\ref{p.tile23}). 

The basic auxiliary fact on 2-attractors is they are uniquely defined, up to affine similarity, by the spectrum of the dilation matrix~$M$ and do not actually depend on the digits~$D$. 
This was first proved in~\cite{Gel}. In Section~3 we give another proof, which makes it possible  
to extend this property to arbitrary attractors whose digits  form an arithmetic progression (Theorem~\ref{th.similar-m}). Thus, the type of 2-attractor depends entirely on  the characteristic polynomial of the matrix~$M$. Moreover, opposite polynomials 
(i.e., the characteristic polynomials of matrices  $M$ and $-M$) generate one and the same 
type. Theorem~\ref{th.uniq} from Section~9 establishes the opposite: 
two attractors corresponding to different and not opposite polynomials 
 are not affinely similar to each other. This fact, whose proof is quite difficult, 
 allows us to find the total number of different  (up to affine similarity) 
 2-attractors in an arbitrary dimension~$d$. This number is equal to the number of different and not opposite to each other integer expanding polynomials of degree~$d$ 
 with the leading coefficient~$1$ and the free coefficient~$\pm 2$. {\em Expanding} 
 means that all roots of the polynomials are larger than one in absolute value.  
 
 This way we obtain the number  $N(d)$ of different 
 2-attractors in the space~$\re^d$: $\ N(2)= 3, \, N(3) = 7, \, N(4) = 21$, etc.
 In Section~10 (Theorem~\ref{th.estimates}) we show that  
 $C_1d^2 \le  N(d) \le C_2 2^{d}$ (the constants are given in the statement of the theorem). 
 The upper bound is derived from the result of Dubickas and Konyagin on the total number of 
 integer polynomials with a given Mahler measure~\cite{DK}. 
 To obtain the lower bound we construct an infinite series of integer expanding polynomials  $p(z) = z^d + a_{d-1}z^{d-1}+ \cdots + a_1 z  \pm 2$ containing a quadratic in $d$
 number of polynomials. Before that, only two series were known in the literature, both of them  are linear in~$d$~\cite{HSV}. In Section~12 we add five new series and a quadratic one among them. We see that the lower and upper bounds on $N(d)$ are far from each other. Numerical results show that at east for  $d\le 8$, our upper bound is tight:   $N(d)$ 
 grows as  $2^d$.  However, this is an open problem to give a proof of this fact or 
 to come up with sharper bounds.

Finally, in Section~11 we analyse the regularity of  2-attractors and compute the 
H\"older exponents in~$L_2(\re^d)$ for the corresponding Haar functions. 
We do it in the  isotropic case for all  $d$  and in the general case for small $d$. 
By the results obtained one can conclude that different 2-attractors always 
have different smoothness. 
 \smallskip 

\smallskip 

\textbf{A brief survey of the literature on 2-attractors.}  
The first examples of 2-attractors were presented in~\cite{G81, B91}
and some properties were studied. The pioneering works~\cite{GH, GM, LW95, LW97}
provided foundations of the theory of tiles and attractors also paid much attention 
to the two-digit case. The topology of plane 2-attractors and combinatorial properties of the corresponding tilings of the plane were studied in~\cite{KL00, KLR, AG, BG} (see 
also Section~6 for more references). In~\cite{Gel} it was shown that there are finitely many 
different, up to cyclic permutations, 2-attractors in~$\re^d$; in~\cite{HSV} 
series of 2-attractors growing linearly in~$d$ were presented. 
All types of 2-attractors in dimensions 2 and 3 (3 and 7 types respectively)
have been found in~\cite{BG}. The monograph~\cite{FG} and papers \cite{GJ, HL, NSVW, Zakh}
study various aspects of plane 2-attractors.  
In~\cite{LW95} six types of plane  2-attractors, up to integer affine similarity, 
were found. This classification was extended in the work~\cite{KL02} to attractors with $3,4$, and $5$ digits. 
In~\cite{Zai} the exponents of regularity of the plane 2-attractors were found. 
Some generalisations for non-integer matrices and related Pisot tilings, 
shift radix systems, etc. were studied in~\cite{KT, ST}. 

Application to wavelets and related issues have been addressed in many papers. 
Here we mention only the monographs~\cite{CHM, Woj, NPS}. We also mention another approach 
to the construction of wavelets from tiles, not involving  MRA, see~\cite{BL, BS, CM, DLS, Mer15, Mer18}.

More references can be found in the corresponding sections.

\smallskip 

\textbf{Notation.} We always assume the basis in  $\re^d$ to be fixed and 
identify linear operators with the corresponding matrices.  
We denote the identity matrix by $I$, the characteristic polynomial of a matrix~$A$
by $p(\lambda) = {\rm det}\, (\lambda I - A)$.  
We use bold symbols for vectors and standard symbols for their components, i.e., 
 $\bx = (x_1, \ldots , x_d)$. By  $|X|$ we denote the Lebesgue measure of a set 
 $X \subset \re^d$.

\bigskip 

\newpage

\begin{center}
\textbf{2. 2-attractors in $\re^2$}
\end{center}
\bigskip 

\begin{figure}[ht!]
\begin{minipage}[h]{0.3\linewidth}
\center{\includegraphics[width=1\linewidth]{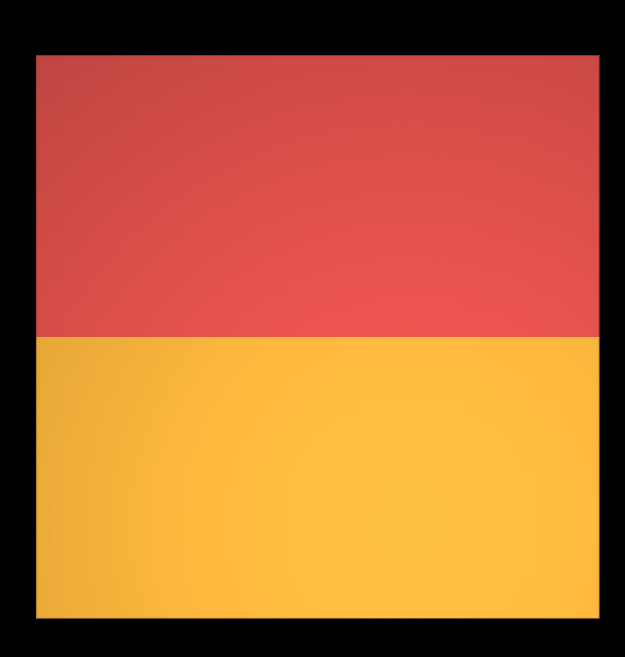}} Square \\
\end{minipage}
\hfill
\begin{minipage}[h]{0.3\linewidth}
\center{\includegraphics[width=1\linewidth]{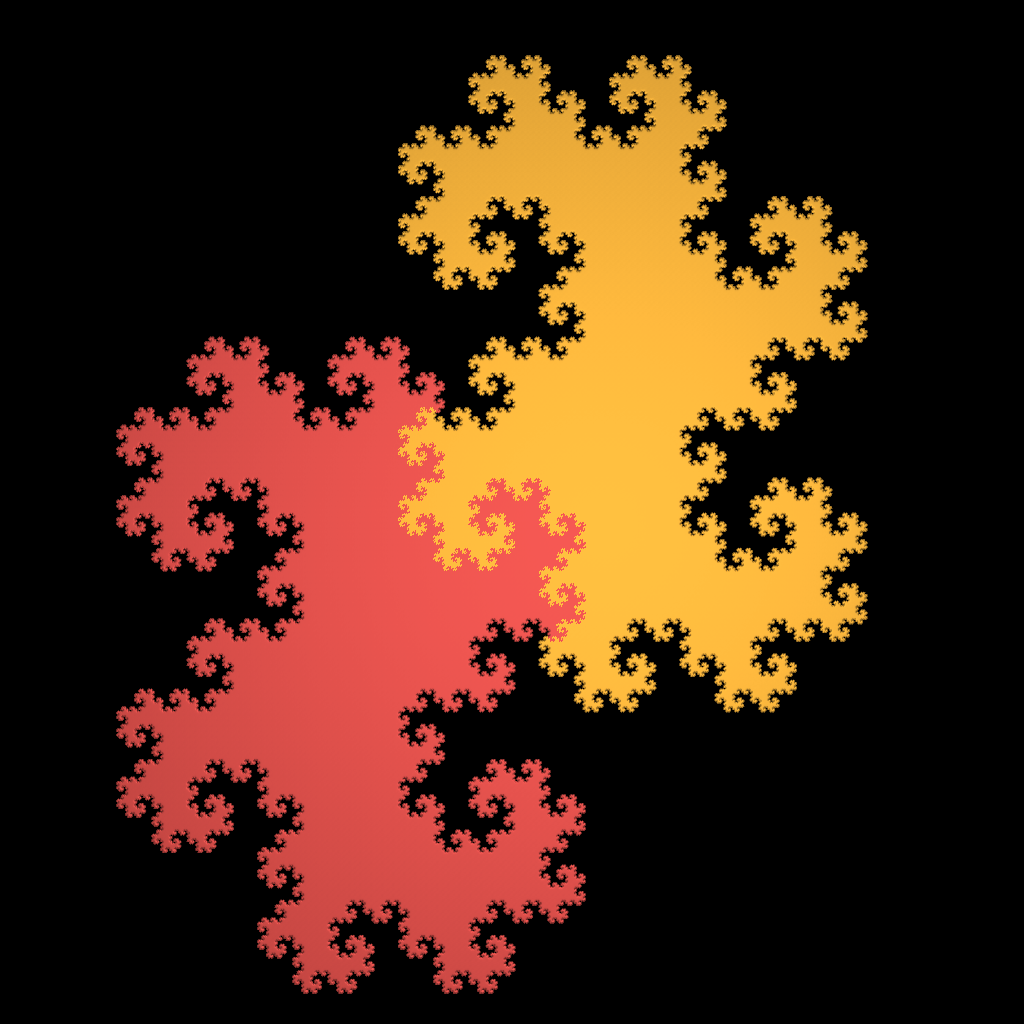}} Dragon \\
\end{minipage}
\hfill
\begin{minipage}[h]{0.3\linewidth}
\center{\includegraphics[width=1\linewidth]{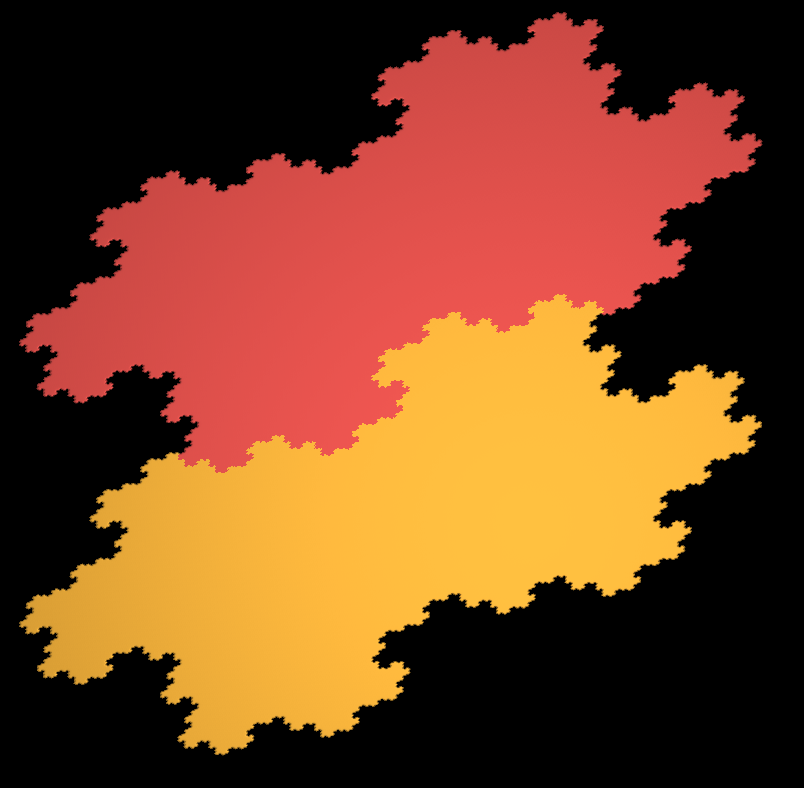}} Bear \\
\end{minipage}
\caption{The partitions of the plane 2-attractors to two affinely-similar parts}
\label{pic1_planar}
\end{figure}

\begin{figure}[ht!]
\begin{minipage}[h]{0.3\linewidth}
\center{\includegraphics[width=1\linewidth]{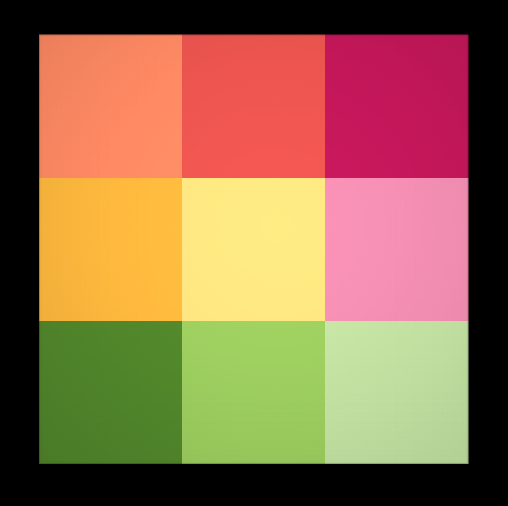}} Square \\
\end{minipage}
\hfill
\begin{minipage}[h]{0.3\linewidth}
\center{\includegraphics[width=1\linewidth]{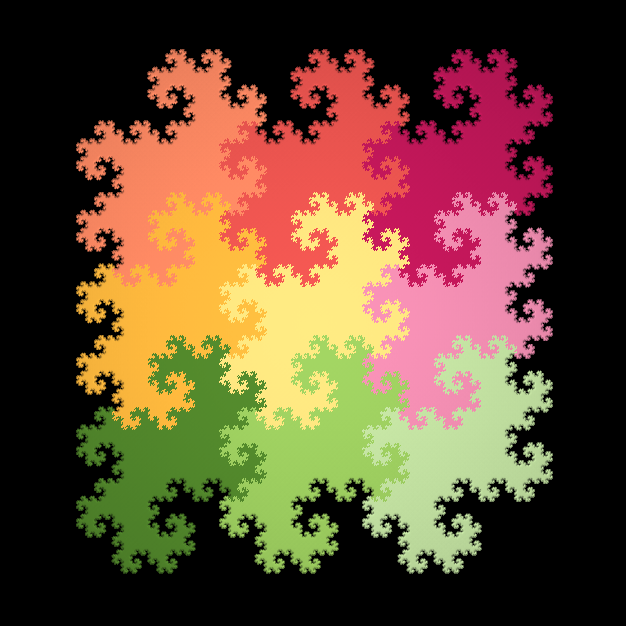}} Dragon \\
\end{minipage}
\hfill
\begin{minipage}[h]{0.3\linewidth}
\center{\includegraphics[width=1\linewidth]{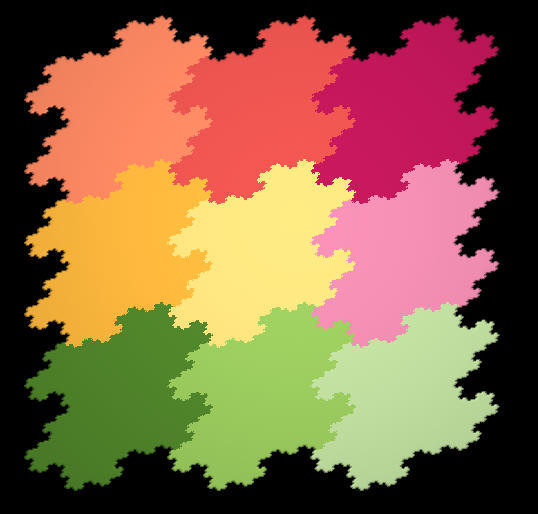}} Bear \\
\end{minipage}
\caption{Tiling by plane  2-attractors $\re^2$}
\label{pic2_planar}
\end{figure}

We begin with the case $d=2$. On the two-dimensional plane there are exactly three, up to 
affine similarity, 2-attractors~\cite{BG, LW95, Zai}. 
We call them  \textit{square}, \textit{dragon} 
 (some papers also use the term ``twindragon'') and \textit{bear} 
 (``tame twindragon''). They are all tiles and all homeomorphic 
 to the disc~\cite{BW}, but not metrically equivalent~\cite{P20}. 
 Their partitions to two parts affinely similar to them are shown 
 in Fig. \ref{pic1_planar}. The tilings by their integer shifts are shown in Fig. \ref{pic2_planar}. The construction of Haar functions and wavelet systems by the 
 plane 2-attractors was realized in~\cite{LW95, Woj, GM}, the smoothness of those Haar functions was analysed in~\cite{Zai}.

The first systematic study of the plane  2-attractors was presented in~\cite{G81, B91}. 
In~\cite{LW95} the matrices of plane  2-attractors are classified 
up to integer similarity. In this case matrices $A$ and $B$ are considered to be equivalent
if there exists an unimodular integer matrix $P \in \mathop{GL}_2(\z)$ such that  $P^{-1}AP = B$. There are $6$ such classes that are unified in three classes of affine isomorphic matrices. Also in~\cite{FG} the case of non-integer digits was studied and it was found for which digits 
there is a tiling of the plane in each of those six cases. 

In case of three and more digits the classification problem is much more difficult because an attractor is already not uniquely defined by the spectrum of the dilation matrix. 
Nevertheless, for small determinants $m = 3, 4, 5$, the work~\cite{KL02} describes all classes of integer affine similarity (but not arbitrary affine similarity) of matrices for each of possible characteristic polynomials.

\bigskip

\begin{center}
\textbf{3. The spectrum of the matrix defines the 2-attractor}
\end{center}
\bigskip 

One of remarkable properties of 2-attractors is that its geometry does not depend on the set of digits~$D$. This was first observed in~\cite{Gel}.  
In Theorem~\ref{th.similar-2} we show that a 2-attractor is uniquely defined, up to an affine similarity, by the spectrum of the dilation matrix. Such a uniqueness cannot be extended to 
a larger number of digits, say, to 3-attractors. Nevertheless, for special digit sets, 
for example, for progressions, a generalization is possible (Section~4). 

We begin with the observation that one digit in~$D$ can be assumed to be equal to zero. 
Indeed,  $\sum_{i=1}^{\infty} M^{-k}\bd_k \, = \, 
\bigl( \sum_{i=1}^{\infty} M^{-k}\bigr)\, \bd_0 \, + \,    \sum_{i=1}^{\infty} M^{-k}(\bd_k - \bd_0)$. Therefore,  if one replaces $\{\bd_0, \bd_1\}$
by $\{\nul, \bd_1 - \bd_0\}$, then the attractor $G$ is translated by vector
$\bigl( \sum_{i=1}^{\infty} M^{-k}\bigr)\, \bd_0$. Thus, in the sequel we assume that $
\bd_0 = 0$, and now everything depends on the choice  $\bd_1$.

We need the following  simple auxiliary fact:  
\begin{lemma}\label{l.commut}
For an arbitrary $d\times d$ matrix $M$ and for arbitrary points  $\ba, \bb \in \re^d$ 
such that  $\ba$ does not belong to an eigenspace of~$M$, there exists a matrix~$C$
commuting with  $M$ and taking  $\ba$ to $\bb$. 
\end{lemma}
{\tt Proof.}  
We pass to the Jordan basis of the matrix 
$M$ and consider one Jordan block~$\Lambda$ of size $s$ corresponding to an 
eigenvalue~$\lambda$. Denote by $\ba'$ and $\bb'$ the $s$-dimensional components of the 
vectors  $\ba$ and $\bb$ corresponding to this block. Denote by $\cH$ the set of Hankel upper-triangular matrices of size~$s$, 
i.e., matrices of the form
$$
\left(
\begin{array}{ccccc}
\alpha_1 &  \alpha_2 & \cdots & \alpha_{s-1} & \alpha_{s} \\
0 &  \alpha_1 & \alpha_2 & \cdots & \alpha_{s-1} \\
0 &  0 & \alpha_1 & \alpha_2 & \cdots \\
\vdots & \vdots & \ddots & \ddots& \vdots \\
0 & 0  & \cdots & 0  & \alpha_1
\end{array}
\right)\, . 
$$
Note that the set $\cH$ is a linear space and is a multiplicative group. 
All elements of $\cH$ commute with $\Lambda$. The set $U \, = \, \{X\ba' \ | \ X \in \cH\}$ is a linear subspace of  
$\re^s$ invariant with respect to every matrix from~$\cH$. 
Let us prove that $U =\re^s$. If this is not the case, then~$U$ 
is an invariant subspace for each matrix from~$\cH$, in particular, for the matrix~$\Lambda$. 
 Therefore, $U$ coincides with one of spaces  
 $U_j = \bigl\{(x_1, \ldots , x_j, 0, \ldots , 0)^T \in \re^s\bigr\}$
 (they form the complete set of  invariant subspaces of~$\Lambda$). 
 Hence, $\ba' \in U_j$ for some~$j\le s-1$. Consequently, the vector $\ba$
 belongs to an eigenspace of the matrix~$M$ generated by all vectors of its Jordan basis 
 except for the last $s-j$ vectors, which correspond to the block~$\Lambda$. 
 This contradicts to the condition on~$\ba$. Thus, $U =\re^s$. 
 Hence, there exists a matrix $X \in \cH$ for which $X\ba' = \bb'$. 
  Let us remember that $X$ commutes with~$\Lambda$. 
  Having found such a matrix~$X$ for each Jordan block of the matrix~$M$
  we compose a block matrix~$C$ for which  $CM  = MC$ and $C\ba = \bb$.

{\hfill $\Box$}
\smallskip

\begin{prop}\label{p.similar-2}
All 2-attractors generated by one matrix~$M$ are affinely similar. 
\end{prop}
{\tt Proof.}  Let a 2-attractor~$G$ be generated by a matrix $M$ and digits $\{\nul, \ba\}$, 
and a 2-attractor $F$ be generated by the same matrix and digits $\{\nul, \bb\}$. 
By the definition of attractors,  
$\be$ does not lie in an eigenspace of~$M$. Therefore, Lemma~\ref{l.commut} 
gives a matrix~$C$ such that  $CM = MC$ and $C\ba = \bb$. 
Every point of $F$ has the form  
$$
\by \ = \  \sum_{j}M^{-k_j}\bb \ = \ \sum_{j}M^{-k_j}C\ba \ =  \ 
\sum_{j}CM^{-k_j}\ba \ = \ C\sum_{j}M^{-k_j}\ba
$$ 
Thus, $\by = C\bx$, where $\bx\in G$. Hence,  $F  = CG$. 
 
{\hfill $\Box$}
\smallskip

\begin{remark}\label{r.8}
{\em According to Proposition~\ref{p.similar-2}, up to an affine similarity, 
a 2-attractor depends only on the expanding matrix~$M$ 
and does not depend on the digit set~$D$. 
As we agreed, one of the digits in~$D$ is equal to~$\nul$. The second digit can be an arbitrary 
integer point out of eigenspaces of~$M$. Moreover, this point can be non-integer, 
although this contradicts to the definition of an attractor.  Anyway the generated attractor~$G$ is affinely similar to a 2-attractor with integer digits. 
To see this it suffices to take an arbitrary point $\bb \in \re^d$ in the proof 
of Proposition~\ref{p.similar-2}, nothing changes. That is why, in what follows, we will sometimes allow the digit~$\bd_1$ to be non-integer in proofs. 

If the converse is not stated, assume that~$\bd_1$
is the first basis vector $\be \, = \, \be_1 \, = \, (1,0, \ldots ,0)^T   \in \re^d$.  

}
\end{remark}

It is known~\cite[Proposition 2.5]{KL00} that an integer expanding matrix with a prime 
determinant cannot have multiple eigenvalues. 
Hence, every such matrix is similar to a diagonal matrix. Therefore, two such matrices 
with the same spectrum are similar. By applying Proposition~\ref{p.similar-2} we obtain  

\begin{theorem}\label{th.similar-2}
A two-digit attractor is uniquely defined, up to an affine similarity,  by the spectrum of the matrix~$M$. 
\end{theorem}
{\tt Proof.}  Let $M= Q^{-1}\Delta Q$, where  
$\Delta$ is a diagonal matrix. Then for  $\bx \in G$, we have 
$$
\bx \ = \ \sum_{k\in \n}M^{-k}\ba_k \ = \ 
\sum_{k\in \n}Q\Delta^{-k}Q^{-1}\ba_k \ = \ 
Q\, \sum_{k\in \n}\Delta^{-k}\, (Q^{-1}\ba_k)\, , 
$$
 where  $\ba_k \in \{\nul, \be \}$ for each  $k$. We see that an attractor with 
 the matrix $M$ and digits $\{\nul, \be\}$ is affinely similar to an attractor 
 with the  matrix $\Delta$ and digits $\{\nul, Q^{-1}\be\}$, 
and  hence, by  Proposition~\ref{p.similar-2}, is affinely similar to an attractor 
with the matrix $\Delta$ and digits $\{\nul, \be\}$. 
 
{\hfill $\Box$}
\smallskip

Thus, a 2-attractor depends on the characteristic polynomial of the matrix~$M$ only. 
The converse is also true: 

\begin{cor}\label{c.similar-2} 
Every integer expanding polynomial with the leading coefficient equal to one and the free 
term   $\pm 2$ generates a unique, up to an affine similarity, 2-attractor. 
\end{cor}
{\tt Proof.}  If $\ p(z) \, = \, z^d + a_{d-1}z^{d-1} + \ldots + a_0$, 
where $a_0 = \pm 2$, is an integer expanding polynomial, then its
 {\em companion matrix}  
\begin{equation}\label{eq.accomp}
M \quad = \quad \left(
\begin{array}{ccccc}
0 &  0 & \cdots & 0 & - a_0\\
1 &  0 & \cdots & 0 & -a_{1} \\
0 &  1 & \cdots & 0 & -a_{2} \\
\vdots & \vdots & \ddots & \ddots& \vdots \\
0 & 0  & \cdots & 1  & -a_{d-1}
\end{array}
\right) 
\end{equation}
has the polynomial~$p$ as a characteristic polynomial. 
Therefore, an attractor generated by the matrix~$M$ with digits $\{\nul, \be\}$, 
corresponds to the polynomial~$p$. By applying Theorem~\ref{th.similar-2} we prove the uniqueness.

{\hfill $\Box$}
\smallskip

\begin{center}
\textbf{4. Attractors generated by arithmetic progressions}
\end{center}
\bigskip 

To what extent can the results of the previous sections be 
generalized to an arbitrary number of digits~$m$? If the digit set~$D$ is arbitrary, then 
any straightforward generalization is impossible already in~$\re^1$. For example, 
for $M=3$, in which case  $m=3$, the attractor generated by the digits $\{0,1,2\}$ 
is a segment while that generated by  $\{0,1,5\}$ is a disconnected fractal set~\cite{P20}. 
So, attractors defined by one matrix but different digits are not necessarily similar. 
The reason is that two-digit sets are always affinely similar to each other  but 
three-digit sets are not. Nevertheless, generalizations are possible 
for special digit sets. For instance, when the digits form an arithmetic progression:
$D \, = \, 
\{\bd_0, \bd_0+\bq, \ldots , \bd_0 + (m-1)\bq\}$, 
where $\bd_0, \bq \in \re^d, \, \bq \ne \nul$. 
Since  the progression  $D \, = \, 
\{\nul, \bq, \ldots , (m-1)\bq\}$ defines the same attractor but shifted by the vector 
$\sum \limits_{k\in \n}M^{-k}\bd_0$, we always assume $\bd_0=\nul$.

 \begin{theorem}\label{th.similar-m}
All attractors generated by a given matrix~$M$ are affinely similar provided the digits 
form an arithmetic progression.   
\end{theorem}
{\tt Proof.}  We assume  $\bd_0 = \nul$. 
If the matrix $C$ commuting with~$M$ (Lemma~\ref{l.commut})
maps the vector $\bq_1$ to $\bq_2$, then it maps the progression ~$D_1 = \{\nul, \bq_1, \ldots \, (m-1)\bq_1\}$ to the progression~$D_2 = \{\nul, \bq_2, \ldots \, (m-1)\bq_2\}$. 
Then we argue as in the proof of Proposition~\ref{p.similar-2}. 
 
{\hfill $\Box$}
\smallskip
 
Theorem~\ref{th.similar-2}, which states that a 2-attractor is uniquely defined by 
the spectrum of the dilation matrix, is not generalized to an arbitrary number of digits, even to progressions. This is because of multiple eigenvalues that can occur. 
The matrix~$M$ can have nontrivial Jordan blocks and hence, can be defined by its 
spectrum in a non-unique way. One can claim that an integer expanding polynomial has no multiple roots only in the case when its free term is a prime number, see~\cite{KL00}. 
Therefore, the following weakened version of Theorem~\ref{th.similar-2} holds 
for arbitrary {\em prime}~$m$:

\begin{cor}\label{c.similar-m}
If~$p$ is an integer expanding polynomial with a prime free term, then all attractors 
generated by matrices with the characteristic polynomial~$p$ and with digits forming an arithmetic progression are affinely similar. 
\end{cor}
{\tt Proof.}   
Arguing as in the proof of Theorem~~\ref{th.similar-2} we establish an affine similarity of attractors generated by the matrices 
$M$ and $\Delta$. We use the fact that the operator   $Q^{-1}$ 
maps a progression to a progression.  
Then we invoke Theorem~\ref{th.similar-2}.  
 
{\hfill $\Box$}
\smallskip

\begin{center}
\textbf{5. Classification of isotropic  2-attractors}
\end{center}
\bigskip  

We come to one of the main result of our paper: to a complete list 
of isotropic 2-attractors. This list turns out to be surprisingly simple.  
Recall that a matrix is called isotropic if it is similar to an orthonormal matrix multiplied by a number. An attractor with an isotropic dilation matrix is also called isotropic. 

Isotropic attractors and tiles have a special place in wavelets theory and in multivariate approximation theory. On the one hand,  an isotropic 
dilation is the most natural for most of problems. In this case the space is 
expanded ``equally in all directions''.   On the other hand, isotropic dilations are 
much simpler to study. It is isotropic wavelets for which
there are direct generalizations of univariate constructions to the multivariate case.
 This concerns, for example, methods of computation of regularity exponents for
 scaling and wavelet functions, methods of estimating the rate of wavelet approximation, etc.   
That is why most of the literature on multivariate wavelets and subdivisions deal with 
isotropic dilation matrices, see,~\cite{Woj, CHM, KPS} and detailed surveys in those works.  

In this section we find all isotropic  2-attractors. It turns out that 
their variety is not rich: in odd dimensions they are parallelepipeds only, 
in even dimensions direct products of dragons and direct products of bears are added. 
Conclusions from this classification will be drawn later. 

An algebraic polynomial is   referred to  as {\em isotropic} 
if all its roots have the same absolute value. An integer polynomial that has the leading 
coefficient equal to one and a prime free term does not have multiple roots~\cite{KL00}. 
Therefore, an expanding integer matrix with a prime determinant is isotropic if and only if its characteristic polynomial is isotropic. 

 \begin{theorem}\label{th.isotr-odd}
If  $d$ is odd, then all isotropic 2-attractors in~$\re^d$ are parallelepipeds.  
\end{theorem}
{\tt Proof.} If  $d$ is odd, then the matrix~$M$ has at least one real eigenvalue, 
which, by the isotropy, takes a value of  either  $2^{1/d}$ or $-2^{1/d}$. 
Substituting to the characteristic polynomial  $p(z) = z^d \, + \, a_{d-1}z^{d-1}\, + \cdots \, + a_1 z\, + \, a_0 $  we obtain
 $- (z^d + a_0) \, = \, a_{d-1}z^{d-1}\, + \cdots \, + a_1 z$. 
Observe that the number $z^d + a_0$ is integer, 
hence, the number $z = \pm 2^{1/d}$ is a root of an integer polynomial of degree~$\le d-1$. 
The latter is possible only if that polynomial is zero. 
Consequently, $p(z) \, = \, z^d + a_0$. Consider the case  $a_0 = - 2$, 
the case $a_0 = 2$ is completely similar. Since the characteristic polynomial annihilates the matrix~$M$, we have $p(M)\, = \, M^d - 2I = 0$. If $\{\nul, \be\}$ are digits generating an 
attractor $G$, then $M^d\be = 2\be$. Let  
$\ba_k = M^{k-1}\be\, , \  k = 1, \ldots , d$, and let 
$P$ be the parallelepiped with a vertex  $\nul$ with edges~$\ba_1, \ldots , \ba_d$. 
Then $P = M^{-1}P\sqcup M^{-1}(P+\be)$. Hence, the characteristic function of the set  $P$ satisfies the same refinement equation:  
$\varphi(\bx) \, = \, \varphi(M\bx) +  \varphi(M\bx-\be)$, as the characteristic 
function of~$G$. Since any refinement equation has a unique, up to multiplication by a constant,  solution, we conclude that $\chi_{P} = \mu \chi_{G}$, where  $\mu \in \re$
is a constant. Since both those functions takes the only values  $0$ and $1$, we obtain
$\mu=1$ and $\chi_{P} =  \chi_{G}$.

{\hfill $\Box$}
\smallskip

Thus, in odd dimensions there are no isotropic 2-attractors except for 
parallelepipeds. In even dimensions, we have a different situation. 
Already in~$\re^2$ two new 2-attractors appear: the dragon and the bear. 
It turns out that there are the same three types in each even dimension.

 \begin{theorem}\label{th.isotr-even}
If $d = 2k$ is even, then in~$\re^d$ there are, up to an affine similarity, 
exactly three isotropic 2-attractors: a parallelepiped, a direct product of $k$ dragons, and a 
direct product of $k$ bears.  
\end{theorem}
This means that in a suitable basis in~$\re^d$, the isotropic 2-attractor has the form~$G = (K, \ldots , K)$, where $K$ is either a rectangle
or a dragon or a bear. The $i$th component $K$
takes the positions $x_{2i}, x_{2i+1}$. 
The crucial point in the proof is the following  
\begin{lemma}\label{l.reduction}
Every integer isotropic polynomial of degree~$d = 2k$ 
with the leading coefficient equal to one and with a prime free coefficient is equal to 
$q(z^k)$, where  $q(t) \, = \, t^{2} + at + b$ is a quadratic integer isotropic polynomial. 
\end{lemma}
{\tt Proof.}  Let the polynomial have the form $z^{2k} + a_{2k - 1}z^{2k - 1} + \ldots + a_1z \pm p = 0$, where  $p$ is a prime number and $a_1, \ldots, a_{2k -1} \in \z$. 
The isotropy implies that the absolute values of all roots are equal to $p^{1/d}$. 
If there is at least one real root among them, then it is $\pm p^{1/d}$ and by repeating  the   proof of Theorem~\ref{th.isotr-odd} we conclude that the polynomial is equal to 
 $z^{2k} \pm p$. The corresponding quadratic integer isotropic polynomial is 
  $q(t) = t^2 \pm p$, which completes the proof. Now assume that all the roots are not real.  
  Then they are divided to pairs of conjugates: $z_1 = \overbar{z_2}, z_3 = \overbar{z_4}, \ldots,  z_{2k-1} = \overbar{z_{2k}}$. The product of roots inside one pair is $z_{2m - 1}z_{2m} = |z_{2m - 1}|^2 = p^{1/k}$, $m = 1, \ldots, k$. 
Let us show that $a_{2k - m} = a_{m} = 0$ for $m = 1, \ldots, (k - 1)$. This will imply that  the polynomial has the form  $z^{2k} + a_k z^{k} \pm p$. By changing the variables  $t = z^{k}$ the polynomial becomes  $q(t) = t^2 + a_k t \pm p$. It remains isotropic since all the roots of  $k$th degree still have the same modulus. The proof will be completed. 

To establish the equality  $a_{2k - m} = a_{m} = 0$ for $m = 1, \ldots, (k - 1)$, 
we first show that  $a_{m} = p^{(k - m) / k} a_{2k - m}$ for $m = 1, \ldots, (k - 1)$. 
Then the numbers  $a_{2k - m}, a_m$ are both zeros since they are integer and the number  $p^{(k - m) / k}$ is irrational for each $m = 1, \ldots, (k - 1)$. 

By the Vieta formulas we have  
\begin{equation}\label{viet}
\begin{aligned}(-1)^m a_m =  \sum \limits_{i_1 < i_2 < \ldots < i_{2k - m}} &{ z_{i_1} \ldots z_{i_{2k - m}}} \\
 (-1)^m a_{2k - m} =  \sum \limits_{j_1 < j_2 < \ldots < j_{m}} &{z_{j_1} \ldots z_{j_{m}}}
 \end{aligned}
 \end{equation}
To each term in the second sum $s = z_{j_1} \ldots z_{j_{m}}$ we associate 
a term $s_1$ in the first sum as follows. First we set $s_1$ to be the complement of~$s$ that consists of those $z_i$ not presented in~$s$. 
Then for those  $z_j$ which do not have their conjugate pairs in $s$ 
we replace in $s_1$ their conjugate pairs back to $z_j$. 
For example, if $k = 4, m = 3$, then the term $z_1z_2z_5$ is associated to  $z_3 z_4 z_5 z_7 z_8$. Let us show that $s_1 = s p^{(k - m) / k}$. Then by formulas  (\ref{viet}) we have  
$a_{m} = p^{(k - m) / k} a_{2k - m}$, which concludes the proof.  
The difference between terms $s_1$ and $s$ is only in those pairs of conjugate rooots 
that are either both presented in~$s$ or both not presented in~$s$. 
Let there be $a$ and  $b$ such pairs respectively. All products of conjugate roots are equal to  $p^{1/k}$, hence, $\frac{s_1}{s} = (p^{1/k})^{(b - a)}$. 
There are in total $m$  multipliers presented in~$s$, those are the $2a$ complete pairs and 
$k - a - b$ incomplete ones. This yields that  $m = k + a - b$, $b - a = k - m$, 
and the desired formula $s_1 = s p^{(k - m) / k}$ is proved.

{\hfill $\Box$}
\smallskip

{\tt Proof of Theorem~\ref{th.isotr-even}}. If  $p$ is an integer isotropic polynomial of degree $d = 2k$, then by Lemma~\ref{l.reduction} we have  
$p(z) = q(z^k)$, where  $q(t) \, = \, t^{2} + at + b$ is a quadratic integer isotropic polynomial. Let 
$$
A \ = \ 
 \left(
\begin{array}{cc}
0 &   - b\\
1 &   -a 
\end{array}
\right) 
$$
be the companion matrix of the polynomial~$q$. Since~$|{\rm det} A | = |b| = 2$ and 
the polynomial~$q$ is isotropic, then the moduli of its roots are~$\sqrt{2}$, 
hence, the matrix~$A$
is expanding. It is proved in~\cite{BG} that every two-digit attractor in~$\re^2$
is either a rectangle or a dragon or a bear.   
Consequently, the matrix~$A$ generates one of those attractors. 
Consider the following matrix: 
\begin{equation}\label{eq.d-even}
M \quad = \quad \left(
\begin{array}{ccccc}
0 &  I & 0  & \cdots  & 0\\
0 &  0 &  I & \cdots  & 0 \\
\vdots & \vdots & \ddots & \ddots& \vdots \\
0 &  0 & \cdots & 0 & I \\
A & 0  & \cdots & 0  & 0
\end{array}
\right) \,  
\end{equation}
defined by the $2\times 2$ blocks: each zero denotes the zero  
$2\times 2$  matrix, $I$ is the identity $2\times 2$ matrix. 
Thus, $M$ is a $d\times d$ matrix and formula~(\ref{eq.d-even})
is defined by $k^2$ blocks of size $2\times 2$. We present an arbitrary vector  $\bx \in \re^d$
as a family of $k$ vectors of dimension~2, i.e., 
$\bx = (x_1, \ldots , x_k), \ x_i \in \re^2$. Then the equation to an eigenvector~$M\bx = \lambda \bx$
becomes the system of equations~$Ax_1 = \lambda x_k, \, 
x_{i+1} \, = \, \lambda x_{i}, \, 
i = 1, \ldots , k-1$. Therefore, $Ax_k \, = \, \lambda^kx_k$, 
and hence,  $\lambda^k$ is an eigenvalue of~$A$. Thus, 
$q(\lambda^k) = 0$, which yields  $p(\lambda) = 0$. 
Thus, the set of roots of the polynomial $p$ (there are exactly $d$ roots, all are simple!) 
coincides with the set of eigenvalues of~$M$.  Therefore, $M$ has a characteristic 
polynomial~$p$. Now we use Theorem~\ref{th.similar-2} and conclude that it is enough 
to prove that the attractor generated by the matrix~$M$, 
is a product of~$k$ attractors of the same type  (rectangle, 
dragon, or bear). Then it will remain to remark that a product of several 
rectangles is a parallelepiped. 

Denote by  $L_i$ the two-dimensional subspace corresponding to the $i$th  
block of the matrix~(\ref{eq.d-even}), $i=1, \ldots , k$. 
Each vector $\bx \in \re^d$ is presented as a sum
$\bx = \sum_{i=1}^k x_i, \, x_i \in L_i$, which will be denoted as~$(x_1, \ldots , x_k)$. 
By $(X_1, \ldots , X_k) = \sum_{i=1}^k X_i$ we denote the direct sum of the sets $X_i \subset L_i, \, i = 1, \ldots , k$, 
which consists of vectors $(x_1, \ldots , x_k), \, x_i \in X_i, \, i = 1, \ldots , k$. 

Let a matrix  $A$  and digits $\{0,a\}$ generate an attractor~$K$
on the plane. Let us show that the matrix $M$ and the digits $\{\nul , \ba\}$, 
where $\ba = (0, \ldots , 0, a)$ ($k$ two-dimensional blocks) generate the attractor~$G = (K, \ldots , K)$. Since  $K$ is  a rectangle, or dragon, or bear, the proof will be completed. 
We have 
\begin{equation}\label{eq.d-even-inverse}
M^{-1} \quad = \quad \left(
\begin{array}{ccccc}
0 &  0 & 0  & \cdots  & A^{-1}\\
I &  0 &  0 & \cdots  & 0 \\
0 &  I & \ddots & 0 & 0 \\
\vdots & \vdots & \ddots & \ddots& \vdots \\
0 & 0  & \cdots & I  & 0
\end{array}
\right) \, . 
\end{equation}
 Consequently, $M^{-1}G \, = \, (A^{-1}K, K, \ldots , K)$ and  
 $$
 M^{-1}(G+\ba)\ = \  M^{-1}G \, + \, M^{-1}\ba\ = \
  (A^{-1}(K+a), K, \ldots , K).
 $$ 
Since  $ A^{-1}K \, \sqcup \, A^{-1}(K+a)\, = \, K$, we see that 
$M^{-1}G \, \sqcup \, M^{-1}(G+\ba) \, = \, G$. 
Therefore,  $G$ is an attractor generated by the matrix~$M$ and by the digits~$\{\nul, \ba\}$. 
Replacing~$\ba$ with any other digit we obtain an affinely similar attractor, which concludes the proof.

{\hfill $\Box$}
\smallskip

\begin{remark}\label{r.diff}{\em 
Thus, in an even dimension there are only three types of isotropic 2-attractors.  
Formally we have not proved yet that those types are not affinely similar. 
This will follow from results of Section~6, where we show that they are not only  
affinely similar but not  $C^1$-diffeomorphic. 
}
\end{remark}
\medskip 

The proofs of Theorems~\ref{th.isotr-odd} and \ref{th.isotr-even}
are constructive. They not only classify all isotropic  2-attractors, but also give the method 
to construct them. If we choose arbitrary bases in the two-dimensional subspaces 
$L_1, \ldots , L_k$, then we obtain the set $(G_1, \ldots , G_k)$, where all $G_i$ 
are arbitrary dragons independent of each other. 
For every dragon  $G_1, \ldots G_k$, their direct product is a $d$-dimensional isotropic attractor. The same is true for a product of 
$k$ arbitrary bears. Let us stress that we multiply either  $k$ dragons or $k$ bears, 
there are no mixed products. 

Since the choice of the dragons  (or bears) is arbitrary, the set 
of $d$-dimensional isotropic 2-attractors has some diversity. 
However, it is achieved due to the choice of an arbitrary basis. In each even dimension, there are only three isotropic 2-attractors up to an affine similarity. This fact leads to the conclusion that  it is necessary to study anisotropic attractors as well. 
As for today, the literature on their applications to wavelets and to approximation algorithms 
is not very extensive, we can mention~\cite{Bow, CHM, CGRS, CGV,  CM, CP}. 
\bigskip

\begin{center}
\textbf{6. The topology of 2-attractors}
\end{center}
\bigskip 

Even basic topological properties of self-affine attractors are 
quite difficult to analyse. Many works~\cite{KL00, KLR, AG, BG, GH, HSV, DL} 
studied the problem of connectedness of attractors. In~\cite{KL00, KLR, AG} 
it was shown that all attractors in dimensions   $d = 2, 3, 4$, with digit sets 
forming an arithmetic progression are connected. The boundaries of attractors and the structure of their neighbours were studied in~\cite{CT, AT, B10, TZ, AL, ALT}.

In~\cite{HSV} it was shown that the plane 2-attractors (the dragon and the bear) are connected.  Moreover, in~\cite{BW} it was proved that they are both homeomorphic to a disc. 
Since this is obvious for a rectangle, we obtain that 
 {\em each plane 2-attractor is homeomorphic to a disc}. 
By using our classification of isotropic 2-attractors  (Theorems~\ref{th.isotr-odd}
and~\ref{th.isotr-even}) and keeping in mind that the product of  $k$ compact sets 
homeomorphic to a disc is homeomorphic to a $2k$-dimensional ball we obtain 
\begin{cor}\label{c.gom}
For each~$d\ge 2$,  all  $d$-dimensional isotropic 
2-attractors are homeomorphic to a ball. 
\end{cor}
For anisotropic 2-attractors, this statement may fail even in the three-dimensional 
space~\cite{B10}. A general (quite complicated) algorithm to verify if an attractor is homeomorphic  to 
a ball was presented in~\cite{CT}. Attractors with arbitrary many digits may not be 
homeomorphic to the ball, the corresponding examples in~$\re^3$ can be found in~\cite{CT} (Propositions  8.7 and  8.11, the boundary of an attractor is a torus). 
However, in the isotropic case the question arises on a stronger topological equivalence. 
In particular, are isotropic  2-attractors diffeomorphic or at least  
{\em metric equivalent}, i.e., homeomorphic by a bi-Lipschitz mapping? 
In this section we obtain a negative answer by showing that all the three
plane 2-attractors (the rectangle, the dragon, and the bear) 
are not bi-Lipschitz equivalent. Hence, they are not  $C^1$-diffeomorphic. Then Theorem~\ref{th.isotr-even} implies that all isotropic 2-attractors in  $\re^{2k}$ 
are not bi-Lipschitz equivalent to each other (in $\re^{2k+1}$ there are only parallelepipeds,  which are obviously all diffeomorphic to each other).  

\begin{defi}\label{d.lip}
Two compact sets~$G_1, G_2 \subset \re^d$ are called metrically equivalent or 
bi-Lipschitz equivalent if they are homeomorphic by a bi-Lipschitz map 
 $f: \re^d \to \re^d$.
The bi-Lipschitz property means that there exists a constant $C > 0$ such that 
$C^{-1}|\bx_1 - \bx_2| \, \le \, |f(\bx_1) - f(\bx_2)| \, \le \, C\, |\bx_1 - \bx_2|$
for all~$\bx_1, \bx_2 \in \re^d$.  
\end{defi}
Note that the bi-Lipschitz property implies its bijectivity. 
If two sets are $C^1$-diffeomorphic, then they are bi-Lipschitz equivalent, 
but not vice versa, since bi-Lipschitz maps can be non-differentiable. 
Although by Rademacher's theorem~\cite{Rad} a bi-Lipschitz map 
is almost everywhere differentiable.  
We establish the non-equivalence of attractors by a characteristic invariant 
with respect to bi-Lipschitz maps. Such an invariant is the 
{\em surface dimension}~$\sigma(G)$. This is a  
supremum of numbers~$r \ge 0$ such that  $|G_{\varepsilon}| \, - \, |G|\, \le \, 
C_{r}\, \varepsilon^{\, r}$ for all~$\varepsilon > 0$, 
where $C_{r}$ does not depend on~$\varepsilon$. See~\cite{P20} on properties and computation of the surface dimension of attractors. 
 \begin{lemma}\label{l.invar}
The exponent $\sigma(G)$ is invariant with respect to bi-Lipschitz maps of~$\re^d$. 
\end{lemma} 
{\tt Proof.} Clearly,  $|G_{\varepsilon}| \, - \, |G|\, = \, 
|G_{\varepsilon}\setminus G|$. Furthermore,  
$f(G)_{\varepsilon}\setminus f(G) \, \subset \, 
 f(G_{C\varepsilon}\setminus G)$. The map  $f$ increases the Lebesgue measure by at most~$C^d$ times. Indeed, the Lebesgue measure is equal to the Hausdorff measure~\cite{Mor}, 
and for the Hausdorff measure, this property is obvious. Therefore, 
$|f(G_{C\varepsilon}\setminus G)|\, \le \, C^d|G_{C\varepsilon}\setminus G|$. 
If  $\sigma(G) > r$, then $|G_{C\varepsilon}\setminus G| \, \ge \, C_r \, |C\varepsilon|^{r}$. Thus, $|f(G)_{\varepsilon}\setminus f(G)| \, 
\le \, C_r \, |C\varepsilon|^{r}\, \le\, C_rC^{r+d} \, \varepsilon^r$. 
Since this is true for all~$\varepsilon >0$, we have 
$\sigma(f(G)) > r$. Thus, the inequality $\sigma(G) > r$ yields   $\sigma(f(G)) > r$.  
Hence,  $\sigma(f(G)) \ge \sigma(G)$. We have shown that a bi-Lipschitz map does not reduce the 
surface dimension. By applying this statement to the inverse map we obtain $\sigma(f(G)) = f(G)$. 

{\hfill $\Box$}
\smallskip

 Lemma~\ref{l.invar} makes it possible to prove the non-equivalence of sets merely by comparing 
 their surface dimensions. 
In the paper~\cite{P20} it was shown that the surface dimension of an isotropic attractor 
is equal to the H\"older exponent of its characteristic function in the space~$L_2$:  
\begin{equation}\label{eq.Holder}
 \alpha(G) \ = \ \alpha(\varphi)\ =  \ \sup_{\alpha \ge 0} \, \Bigl\{ 
 \bigl\| \varphi(\bx + \bh)\ - \ \varphi(\bx)  \bigr\|_1 \ \le \ C\, |\bh|^{\, \alpha}\, , \ \bh \in \re^d\, \Bigr\}\, ,  
 \end{equation}
where $\varphi = \chi_{G}$. The H\"older exponent of a parallelepiped, as well as of any 
polyhedron, is equal to~$\frac12$. The H\"older exponents of the dragon and of the bear have been computed in~\cite{Zai} by applying methods from~\cite{CP}. We come to the following result:

\begin{theorem}\label{th.isotr-nehomog}
For every even~$d$, all the three types of isotropic 2-attractors are homeomorphic but not bi-Lipschitz equivalent. 
\end{theorem} 
Since for odd~$d$, all isotropic  2-attractors are parallelepipeds, 
they are, of course, diffeomorphic. Hence, Theorem~\ref{th.isotr-nehomog}
is applied for even dimensions only.  In the proof we will need an auxiliary fact established in~\cite{P20}.  For an arbitrary subspace 
$V \subset \re^d$, the {\em $L_2$ H\"older exponent of $\varphi$ along~$V$}
is the value
$$
\alpha_{V}(\varphi)\ =  \ \sup_{\alpha \ge 0} \, \Bigl\{ 
 \bigl\| \varphi(\cdot + \bh)\ - \ \varphi(\cdot)  \bigr\|_1 \ \le \ C\, |\bh|^{\, \alpha}\, , \ \bh \in V\, \Bigr\}\, 
 $$
 (all the shifts are along the subspace~$V$). 
\smallskip 

\noindent \textbf{Lemma A}~\cite{CP}. {\em If $\varphi \in L_2(\re^d)$ is an arbitrary
function and the space~$\re^d$ is presented as a direct sum of subspaces~$V_1, \ldots , V_k$,
then $\alpha(\varphi)\, = \, \min \limits_{j=1, \ldots , k}\alpha_{V_j}(\varphi)$. 
}

\smallskip 

{\tt Proof of Theorem~\ref{th.isotr-nehomog}.}  For $d=2$, 
the H\"older exponents of two-dimensional 2-attractors  
are computed in~\cite{Zai}. For the dragon, this exponent is $0.2382...$, 
for the bear, it is equal to $0.3446...$ (for the square, it is, of course, $0.5$). 
Since those exponents are equal to the corresponding surface dimensions,  
Lemma~\ref{l.invar} implies that those sets are not bi-Lipschitz equivalent. 
For $d=2k$, each  isotropic 2-attractor is equal to a direct product of either $k$
dragons or  $k$ bears, or is a parallelepiped. In the latter case 
$\alpha(G)=0.5$. In case of dragons we denote by  $V_1, \ldots , V_k$
two-dimensional planes that contain those dragons. 
Then $\alpha_{V_j}(G)$ is equal to the regularity of a dragon for every~$j$. 
By Lemma~A, $\alpha_{ V_j}(G)$ is equal to the regularity of the dragon $0.2382...$. 
Similarly, the regularity of the direct product of $k$ bears is equal to the regularity of one bear, which is $0.3446...$. Now referring again to Lemma~\ref{l.invar} we conclude the 
non-equivalency of the sets. 

{\hfill $\Box$}
\smallskip

As for the non-isotropic  2-attractors, we can only formulate the following conjecture: 
\begin{conj}\label{conj.10}
If 2-attractors are not affinely similar, then they are not bi-Lipschitz equivalent. 
\end{conj}

\bigskip

\begin{center}
\textbf{7. Convex hulls of  2-attractors}
\end{center}
\bigskip 

Polyhedral attractors have been studied in the literature~\cite{GM, CT, Zai2}. 
Apart from the theoretical interest, they are attractive for applications because they 
generate wavelets and subdivision schemes of a simple  structure. 
It was proved in~\cite{Zai2} that every attractor-polygon in~$\re^2$ (not necessarily convex)  is a parallelogram. A conjecture was made that this assertion is true in an arbitrary dimension: 
all polyhedral attractors are parallelepipeds. We are going to prove this conjecture 
for  $2$-attractors. 
Actually, we will establish a stronger fact and find all  $2$-attractors whose 
convex hull is a polytope. We prove that a convex hull of an arbitrary 
$2$-attractor is an infinite zonotope, which, in turn, a polytope precisely in two cases: 
for parallelepipeds and for products of dragons.  

Let us remember that a {\em zonotope} in $\re^d$ is a Minkowski sum  of a finitely many segments. Every zonotope is convex and is a projection of the $N$-dimensional cube
to the space  $\re^d$, where $N$ is a number of segments. 
In particular, a zonotope always has a center of symmetry. A {\em zonoid} is a limit of a  sequence of zonotopes in the 
Hausdorff metric. We need a generalization of the notion of zonotope to an infinite number of segments. 
\begin{defi}\label{d.zon}
An {\em infinite zonotope} in~$\re^d$ is a Minkowski sum 
of a countable set of segments whose sum of lengths is finite.  
An infinite zonotope is nondegenerate if it does not lie in a proper subspace. 
\end{defi}
Every infinite zonotope is a zonoid but not vice versa. A nondegenerate infinite zonotope is a 
centrally symmetric convex body. 
\begin{prop}\label{p.zon}
A convex hull of a  $2$-attractor generated by a matrix~$M$ 
and digits 
$\{\nul, \ba\}$ is a nondegenerate infinite zonotope: 
$\bar G\, = \, \sum_{k=1}^{\infty} M^{-k}\bar \ba$, where  $\bar \ba$ is 
the segment~$[\nul, \ba]$. 
\end{prop}
{\tt Proof.}  
Since~$\rho(M^{-1})< 1$, it follows that the sum of lengths of the segments $M^{-k}\bar \ba$ is finite.  Hence,  $\bar G$ is a nondegenerate infinite zonotope.  Clearly, $G \subset \bar G$, because~$\{\nul, \ba\} \, \subset \, \bar \ba$. Let us prove the opposite inclusion. 
Every point $\bx \in \bar G$ is presented as a sum $\sum_{k=1}^{\infty} M^{-k} \bx_k$, 
where each point $\bx_k$ belongs to the segment~$M^{-k} \bar \ba$. 
If for at least one $ j$, the point  $\bx_j$ does not coincide with the end of the segment~$M^{-j} \bar \ba$, 
then it is a half-sum of some points $\bx_j', \bx_j'' \, \in \, M^{-j} \bar \ba $. 
Therefore, $\bx$ is a half-sum of points~$\bx_j' + \sum_{k\ne j}^{\infty} M^{-k} \bx_k$  and~$\bx_j'' + \sum_{k\ne j}^{\infty} M^{-k} \bx_k$. Hence, all extreme points of 
$\bar G$ are among points of the sum $\sum_{k=1}^{\infty} M^{-k} \{\nul, \ba\}$, 
i.e.,  among points of the set~$G$. By the Krein-Milman theorem,
every convex compact set is a closure of the convex hull of its extreme points. 
Then $\bar G$, as a a closure of the convex hull of extreme points of the compact set~$G$, 
is contained in~$G$.  

{\hfill $\Box$}
\smallskip

As the first corollary, we obtain the following curious fact on dragons. 
Most likely, it is known but we could not find a reference and therefore 
include its proof. 
 
\begin{prop}\label{p.zon-dragon}
A convex hull of a dragon is a convex octagon. 
\end{prop}
{\tt Proof.}  
The matrix~$M$ of a dragon defines a rotation by $\frac{\pi}{4}$ with multiplication by $\sqrt{2}$. Hence, 
all the segments  $M^{-k}\bar \be$ are located on four straight lines, 
which are the coordinate axes and the bisectors of the coordinate corners. All the segments 
on one line are summed up to one segment. So, the sum 
 $\sum_{k=1}^{\infty}M^{-k}\bar \be$ is a sum of four segments, which  is an octagon. 

{\hfill $\Box$}
\smallskip

It turns out that among all 2-attractors, not necessarily isotropic, 
only the parallelepiped and the dragon, or direct products of dragons, 
have simple convex hulls. 
A direct product of  $k$ dragons (one of the three isotropic 2-attractors in  $\re^{2k}$
by Theorem~\ref{th.isotr-even}) is a polytope with $8^k$ vertices. 
All other 2-attractors have convex hulls which are not polytopes. 

\begin{theorem}\label{th.zon-all}
Among all 2-attractors in $\re^d$,  there exist precisely two types that have convex 
hulls which are polytopes. This is a parallelepiped ($2^d$ vertices) and a direct product of 
$\frac{d}{2}$ dragons ($2^{\frac{3d}{2}}$ vertices). 
 
 Every 2-attractor different from a parallelepiped and from a product of dragons, 
has a convex hull that is an infinite zonotope, which has infinitely many extreme points. 
\end{theorem}
In the proof we need the following technical lemma. 
\begin{lemma}\label{l.zon-segments}
If among the segments generating an infinite zonotope in~$\re^d$, 
there are  infinitely many non-parallel, then that zonotope is not a polytope. 
\end{lemma}
{\tt Proof.} Let us first prove the statement in~$\re^2$. 
Having taken sums of all parallel segments we may assume that all segments in the sequence~$\{\bar \ba_k\}_{k=1}^{\infty}$
have different directions. Let us take an arbitrary index~$j$ and consider the sequence of zonotopes~$\bar G_s \, = \, \sum_{k=1}^{s}
\bar \ba_k$ for $s=j, j+1, \ldots $. Each of them is a  $2s$-gon, 
which has one side parallel and equal to the segment~$\bar \ba_j$.  
Since $\bar G_s$ converges to an infinite zonotope~$\bar G$ in the Hausdorff metric, 
it follows that the limit zonotope has 
a segment parallel and not smaller (by length) to the segment~$\bar \ba_j$ 
on its boundary. Hence, the boundary of~$\bar G$ contains segments parallel to all the segments~$\bar \ba_j, \, j\in \n$.  Therefore, $\bar G$ is not a polygon, 
which concludes the proof in~$\re^2$. 
To probe the theorem in~$\re^d$ it suffices to project all the initial vectors 
to some two-dimensional plane so that infinitely many non-parallel vectors remain after this 
projection. By what we proved above, the projection of  $\bar G$ to this plane 
is not a polygon, hence  $\bar G$ is not a polytope.

{\hfill $\Box$}
\smallskip

{\tt Proof of Theorem~\ref{th.zon-all}.}  
Let us first show that a convex hull of the bear is not a polygon. 
Indeed, in a suitable coordinates, the matrix of the bear defines a rotation by the angle ${\rm arctg}\, (\sqrt{7})$ with the expanding by~$\sqrt{2}$. Since the angle is irrational,
the set~$M^{-k}\bar \be\, , \, k \in \n$, does not contain parallel segments. 
It remains to refer to Lemma~\ref{l.zon-segments}. 
This completes the proof for~$\re^2$. 
Consider now arbitrary~$\re^d, \, d\ge 3$. If the attractor~$G$ 
is isotropic, then it is either a parallelepiped or a product of dragons or a product of 
bears. The first two cases are in the assumptions of the theorem. 
In the latter case the convex hull of the attractor is not a 
polytope, since the convex hull of the bear is not a polygon. 

It remains to consider the case when the attractor~$G$ is not isotropic. 
Denote by~$V$ the linear span in~$\re^d$ of eigenvectors corresponding to all 
largest by modulus eigenvalues of the matrix~$M^{-1}$ 
(let us recall that there are no multiple eigenvalues). 
If some eigenvalue is not real, we take the real and imaginary parts of the 
corresponding eigenvector. The vector~$\be$ does not belong to~$V$
because it does not belong to invariant subspaces of the matrix~$M$. 
Hence, the sequence $\{M^{-k}\be\}_{k=1}^{\infty}$ 
approaches the subspace~$V$ as $k \to \infty$ but not reaches it 
since 
$M^{-1}$ is nondegenerate. Therefore, the sequence of segments~$\{M^{-k}\bar \be\}_{k=1}^{\infty}$ contains infinitely many non-parallel segments. Invoking now Lemma~\ref{l.zon-segments} we complete the proof.

{\hfill $\Box$}
\smallskip

Thus, according to Theorem~\ref{th.zon-all},  only parallelepipeds and products of dragons have simple convex hulls. 
\begin{remark}\label{r.bear-drag}
{\em It is interesting that the bear is more regular attractor than the dragon: the 
regularity exponent of the bear in~$L_2$ is equal to $0.3446...$ 
while for the dragon, this is  $0.2382...$. 
Nevertheless, the convex hull of the bear is an infinite zonotope while that 
of the dragon is an octagon.
}
\end{remark}

Now we are able to classify all simple 2-attractors and prove the conjecture from~\cite{Zai2} in case of 2-attractors. 

{\em A polyhedral set}
is a union of a finite number of polyhedra. 
\begin{cor}\label{c.poly}
If a 2-attractor is a polyhedral set, then this is a parallelepiped.  
\end{cor}
{\tt Proof.} By Theorem~\ref{th.zon-all}
if a 2-attractor is a polyhedral set, then this is either a parallelepiped or a product of dragons. The dragon, however, is not a polyhedral set since its H\"older exponent in $L_2$ is strictly less than~$0.5$.

{\hfill $\Box$}
\smallskip

In the proof of Theorem~\ref{th.zon-all} we have shown that a  convex hull of 
a non-isotropic 2-attractor is never a polytope. This immediately implies 
\begin{cor}\label{c.brick}
If a  2-attractor is a parallelepiped, then it is isotropic. 
\end{cor}

\bigskip 

\begin{center}
\textbf{8. Tiles, Haar bases, and MRA}
\end{center}
\bigskip 

The Lebesgue measure of a self-similar attractor is a natural number. 
If this measure is equal to one, then the attractor is called a tile. 
Integer shifts of a tile cover the space~$\re^d$ with intersections of zero measure. 
The characteristic function~$\varphi$ of a tile has 
orthonormal integer shifts. Among all attractors, the tiles are most important in applications. In this case the function~$\varphi$ generates an MRA and a system of wavelets. 
For 2-attractors  (in this case 2-tiles)
$G(M, D)$, where $D=\{\nul, \ba\}$,  
the Haar basis is generated by $M$-contractions and integer shifts of one 
 function 
$\psi(\bx) \, = \, \varphi(M\bx) \, - \, \varphi(M\bx - \ba)$. 
Wavelets systems are also generated by tiles in the frequency domain~\cite{DLS, BL, BS, Mer18, Mer15}. In the approximation theory, the case of tiles is also especially important since 
in this case the subdivision scheme converges in  $L_p$. 
If an attractor is not a tile, then integer shifts of~$\varphi$ are linearly dependent, 
therefore, the function~$\varphi$ does not generate an MRA 
and the corresponding subdivision scheme may diverge~\cite{CDM, KPS}. 
There are several criteria to verify that an attractor generated by a matrix~$M$ and by digits~$D$ is a tile~\cite{LW97, GH}. This depends not only on the matrix but also on the digits. 
For example, if a matrix~$M$ and digits~$D=\{\nul, \ba\}$ generate a 2-tile~$G(M,D)$,  
then the 2-attractor $G(M, D')$ generated by the same matrix and by the digits $D'=\{0, 3\ba\}$
is not a tile since it has the Lebesgue measure $3^d|G(M,D)| \, \ge \, 3^d$. 
In dimensions  $d=2, 3$, for every integer expanding matrix $M$, 
there exists a digit set~$D$ for which  $G(M, D)$ is a tile. 
However, already in~$\re^4$ A.Potiopa in 1997 presented an example of matrix~$M$
for which such a digit set does not exist, i.e., that matrix does not generate a tile~\cite{LW96, LW99}. 
 
 In case of two digits,  the type of 2-attractor depends  neither on  digits nor on 
 the matrix~$M$ but only on the characteristic polynomial of~$M$. 
 That is why it is natural to formulate the problem in a different way:
\smallskip 

\noindent \textbf{Problem.} {\em Does always exist a tile affinely similar to a given 2-attractor~?} 

\smallskip   
 
Since a 2-attractor is uniquely defined, up to an affine similarity, by 
an expanding polynomial, the problem can be formulated in a different way: 
{\em Can every expanding polynomial with   $\, a_d = 1, \, |a_0| = 2$ generate a tile~?}
 \smallskip 

We conjecture that the answer is affirmative and, moreover, 
is attained for the companion matrix~(\ref{eq.accomp}) of that polynomial 
and for digits $D=\{\nul, \be\}$, where $\be = \be_1$ is the first basis vector.  
\begin{conj}\label{conj.20}
For an arbitrary expanding polynomial with the leading coefficient~$1$ 
and a free term $\pm 2$, its companion matrix~$M$ and the digits~$D=\{\nul, \be\}$ 
generate tile. 
\end{conj}
We immediately observe that this conjecture is true in at least two cases. 
The first one is the case of small dimensions. 
 \begin{prop}\label{p.tile23}
Conjecture~\ref{conj.20} is true for dimensions  $d=2,3$. 
\end{prop}
{\tt Proof}. In the dimension $d=2$ there are, up to an equivalence, three polynomials 
generating 2-attractors  (Section~2); in dimension $d=3$ there are seven such polynomials  (Section~11). Applying to each of them the program {\tt checktile}~\cite{Mej} that realizes the criterion from~\cite{GH} we complete the proof.

{\hfill $\Box$}
\smallskip

The second case is more important. Conjecture~\ref{conj.20} is true for isotropic 2-attractors 
in all dimensions. 
 \begin{theorem}\label{th.tile-iso}
Conjecture~\ref{conj.20} is true for isotropic polynomials. 
\end{theorem}
{\tt Proof}. In an odd dimension 
every 2-attractor is a parallelepiped corresponding to the polynomial~$p(z) = z^d  \pm 2$ (Theorem~\ref{th.isotr-odd}). 
Then one verifies directly that the characteristic function of a cube  $[0,1]^d$ satisfies 
the refinement equation $\varphi(\bx) = \varphi(M\bx) + \varphi(M\bx - \be)$, where $M$ is a companion matrix of the polynomial~$z^d  - 2$.  
Hence, for this polynomial, the 2-attractor is a tile~$[0,1]^d$. 
For the polynomial~$z^d  - 2$, the same refinement equation has another solution:  
the characteristic function of a cube  $\bigl[-\frac13, \frac23 \bigr] \times \bigl[ 0, 1\bigr]^{d-1}$. 
This cube apparently is also a tile. 

In the even dimension $d=2k$ Lemma~\ref{l.reduction} yields   
$p(z)= q(z^k)$, where $q(t) \, = \, t^{2} + at \pm 2$ is a quadratic isotropic 
polynomial. A permutation of coordinates 
 $x_i \, \to \, x_{\sigma(i)}$ by the formula  $\, \sigma(i) \, = \, 2(k-i+1);  \  \sigma(i+k) = 2(k-i)+1, \ i = 1, \ldots , k$, maps the companion matrix of the polynomial~$p$ 
 to the matrix~$M$ defined by formula~(\ref{eq.d-even}). 
 Hence, if  $M$ generates a tile, then the companion matrix also does.  
By Theorem~\ref{th.isotr-even},  the matrix~$M$ generates a  2-attractor, which is 
a direct product of equal two-dimensional attractors: dragons, bears, and rectangles.  
Since each of those two-dimensional attractors is a tile, 
their direct product is also a tile.

{\hfill $\Box$}
\smallskip

\bigskip 

\begin{remark}\label{r.tail}
{\em Theorem~\ref{th.tile-iso} guarantees the existence of a 2-tile 
affinely similar to an arbitrary isotropic 2-attractor. 
In the example of Potiopa~\cite{LW99} the matrix $M$ is isotropic as well
and it does not generate a tile, i.e., there is no suitable digit set~$D$ 
such that  $G(M,D)$ is a tile. 
However, Theorem~\ref{th.tile-iso} claims that for the companion matrix of the 
same characteristic polynomial~$p(z) = z^4 + z^2 + 2$, such a digit set exists: this is 
 $D = \{\nul, \be \}$. We obtain a 2-tile $G(M, D)$. 
By Theorem~\ref{th.isotr-even} this tile is a product of two dragons. Thus, in the example of Potiopa it suffices to merely  pass to another integer basis in $\re^d$
to obtain a 2-tile.  
}
\end{remark}
\begin{remark}\label{r.haar}
{\em Despite the fact that an isotropic 2-tile is a direct product of identical 
bivariate tiles, it generates a Haar system which is different from a direct product of 
bivariate Haar functions. In fact, if $\varphi$ is a product of, say, 
$k$ dragons: $\varphi = \varphi_1 \otimes \cdots \otimes \varphi_k$, 
where each function $\varphi_i$ is a characteristic function of a dragon in~$\re^2$, 
then the Haar function in $\re^d$ is $\varphi(M\bx) - \varphi(M\bx - \be)$, 
it is different from the direct product of  $k$ bivariate Haar functions. 

Let us also note that the classification of isotropic  2-attractors in Section~5 
allows us to find the regularity exponents of the corresponding Haar functions. 
For odd $d$, there are only parallelepipeds, their H\"older exponent in~$L_2(\re^d)$
is $0.5$. For even $d$, the list is complemented by products  of dragons and products of bears.
Their regularities in~$L_2(\re^2)$ are equal to the regularity of dragons and of bears 
respectively, this is  $0.2382...$  and $0.3446...$ respectively~\cite{Zai}.

}
\end{remark}

\bigskip 

\begin{center}
\textbf{9. Can one attractor correspond to different matrices?}
\end{center}
\bigskip 

How many different, up to an affine similarity, 2-attractors are there in~$\re^d$? 
In the isotropic case there are either one or three types depending on the evenness of~$d$. 
For non-isotropic matrices, there are much more types, and their total number grows with the dimension. In the next section we will estimate this number. However, before 
address this issue we have to answer one important question: 
 {\em is there a one-to-one correspondence between 2-attractors and the spectra 
 of expanding matrices with the determinant~$\pm 2$ ?}  In Section~3 this 
 correspondence was established in one direction:  by Theorem~\ref{th.similar-2}, 
 each 2-attractor is uniquely defined by the spectrum of its matrix~$M$, 
 i.e., by a polynomial with the leading coefficient~$1$ and the free term~$\pm 2$. 
Hence, we can estimate the number of such polynomials of  degree~$d$. 
Each of them generates a unique, up to an affine similarity, attractor. 
However, is the converse true? Is it true that each attractor is associated to a unique 
polynomial, i.e., a unique, up to similarity, matrix~$M$? 
Is it possible that we find many proper polynomials, but all of them 
define the same type of attractors, say, parallelepiped? This problem turns out to be non-trivial. First, we show, which is simple, that the opposite matrices  $M$ and $-M$ 
define the same 2-attractor. Then, which is more difficult, that different and non-opposite matrices define different (not affinely similar) attractors. 
We need the following known~\cite[Proposition 2.2]{KL00} fact.
\begin{prop}\label{p.symm} 
Every 2-attractor is centrally symmetric. 
\end{prop}
{\tt Proof.} Let $D = \{\nul , \be\}$. 
Denote  $\bc \, = \, \frac12 \, \sum\limits_{j=1}^{\infty} M^{-j} \be$.  
Then $G \, = \, \bc \, + \, \sum \limits_{j=1}^{\infty} \pm \frac12 M^{-j} \be$. 
The set  $\sum \limits_{j=1}^{\infty} \pm \frac12 M^{-j} \be$ is symmetric about the origin, hence $G$ is symmetric about the point~$\bc$. 

{\hfill $\Box$}
\smallskip

\begin{defi}\label{d.prot}
Algebraic polynomials~$p = \sum \limits_{k=0}^{d} p_k t^k$ and 
$q = \sum \limits_{k=0}^{d} q_k t^k$ are called {\em opposite} is 
$q_k = (-1)^{d-k}p_k\, , \ k = 0, \ldots , d$. 
\end{defi}
The leading coefficients of opposite polynomials are equal. 
The Vieta formulas  imply that the roots of the polynomial $q$ 
are roots of the polynomial~$p$ taken with the opposite signs. 
Thus, the roots of $p$ and $q$ are opposite, which justifies our terminology.  
The matrices $M$ and $-M$ have opposite characteristic polynomials.

\begin{prop}\label{p.prot} 
Opposite polynomials generate the same 2-attractor. 
\end{prop}
{\tt Proof.} The characteristic function~$\varphi$  of an attractor $G$
satisfies the refinement equation~$\varphi(\bx) \, = \, 
\varphi(M\bx) + \varphi(M\bx - \be)$. Since $G \, = \, 2\bc \ -\, G$, it follows that  
$\varphi (2\bc - \bx)\, = \, \varphi(\bx)$. Then 
$\varphi(\bx) \, = \,  \varphi(2M\bc - M\bx) + \varphi(2M\bc - \be - M\bx)$.
Therefore, the set  $G$ is also an attractor with the matrix $-M$
and the digits $-2M\bc\, , \, \be\, - \, 2M\bc$. 

{\hfill $\Box$}
\smallskip
 
It turns out that the converse is also true: 
if polynomials generate identical 2-attractors, 
then they are either equal or opposite. This means that 
a 2-attractor defines the dilation matrix~$M$ up to a sign and to a similarity. 

\begin{theorem}\label{th.uniq} 
If matrices~$M_1$ and $M_2$ (with some sets of digits~$D_1, D_2$)  
generate the same, up to an affine similarity, 2-attractor, then either $M_1 = M_2$ or 
 $M_1 = - M_2$.  
\end{theorem}
The proof is based on the following lemma. For a matrix $M$ 
and for a segment  $\bar \ba = [- \ba, \ba ]\subset \re^d, \, \ba \ne 0$, 
consider the following set of segments  $\cT(M, \bar \ba)\, = \, \bigl\{ M^{-k}\bar \ba, \, k \in \n \bigr\}$, where coinciding segments are counted with multiplicity. 
\begin{lemma}\label{l.dyn}(on the recovery of a dynamical system). 
Two expanding $d\times d$ matrices $M_1, M_2$ without multiple eigenvalues are given. 
If for some segments $\bar \ba_1$ and  $\bar \ba_2$ 
($\bar \ba_i$ is not contained in an invariant subspace of $M_i, \, i = 1,2$)
we have 
  $\cT(M_1, \bar \ba_1) \, = \, \cT(M_2, \bar \ba_2)$, then either 
 $M_1 = M_2$, or $M_1 = -M_2$.    
\end{lemma}
{\tt Proof} see Appendix.  In fact, Lemma~\ref{l.dyn} 
guarantees  a possibility of reconstruction of a linear dynamical system by its 
trajectory {\em even if the history is unknown}, i.e., 
the succession of points is unknown. If we know the history, then the system can be identified by each~$d+1$ successive points  in the trajectory. The problem of identification of a system with an unknown history is significantly more difficult. 
\smallskip 

Now we are able to prove Theorem~\ref{th.uniq}. To this end, we 
use a known result from the theory of refinement equations and present the Fourier 
transform~$\widehat \varphi$ of a characteristic functions of an attractor~$G$ 
as an infinite product of trigonometric polynomials. 
By this form we find the set of zeros of the function~$\widehat \varphi$. 
This set is a union of countably many hyperplanes in~$\re^d$. 
To avoid dealing with hyperplanes we pass to their polars. 
The polar is taken with respect to a unit Euclidean sphere centered at $0$. 
The polar of a hyperplane is one point. By applying Lemma~\ref{l.dyn} we show that the obtained set of zeros defines the matrix~$M$ in a unique way. Therefore, the attractor~$G$
uniquely (up to  $\pm$) defines~$M$. 

\smallskip 

{\tt Proof of Theorem~\ref{th.uniq}}. Suppose a 2-attractor~$G$ is generated by a matrix 
$M$ and by digits $\{\nul, \be\}$; then its characteristic function satisfies a refinement equation  
$\varphi( \bx)\, = \varphi(M\bx) \, + \, \varphi(M\bx -\be)$
(see Section~1, equation~\ref{eq.ref}). The solution of a refinement equation is expressed by the following formula
\begin{equation}\label{eq.masks}
\widehat \varphi (\xi)\quad  = \quad \widehat \varphi (0)\prod_{k=1}^{\infty} \bm \bigl( [M^{*}]^{-k}\bxi\bigr)\, ,  
\end{equation}
where the trigonometric polynomial~$\bm(\bxi)\, = \, \frac{1}{2}\Bigl(1 \, + \, e^{-2\pi i (\be, \bxi) } \Bigr)$ is called a  
{\em mask} of the equation (see, for example, \cite{NPS}). If we show that the set of zeros of the function
$\widehat \varphi (\bxi)$ uniquely, up to a sign, defines the matrix~$M$, 
then everything will be proved. Indeed, the function $\widehat \varphi $ 
is uniquely defined by the set~$G$, therefore, the set of zeros of this function is 
also defined by the set~$G$. 
Hence, the matrix~$M$, up to a sign, will also be defined by the set~$G$. 
Since $(\be, \bxi) = \xi_1$ (the first component of the vector~$\bxi$), 
we see that the set of solutions of the equation  $\bm(\bxi)\, = \, 0$
is a union of countably many parallel hyperplanes   $\xi_1 = \frac12 + n, \, n \in \z$. 
The polars of those hyperplanes form a countable set of points  $\cR = \bigl\{ \frac{2}{1+2n}\, \be_1\bigr\}_{n \in \z}$, where $\be_1 = (1,0, \ldots , 0)^T$ is the first basis vector. 
The set of zeros of the function  $\bm \bigl( [M^{*}]^{-k}\bxi\bigr)$ 
is an image of the set of zeros of the function $\bm \bigl(\bxi\bigr)$ under the action 
of the operator~$[M^*]^k$. 
Hence, the polar of the set of zeros of the function  $\bm \bigl( [M^{*}]^{-k}\bxi\bigr)$ 
is  $M^{-k}\cR$. The set of zeros of the product~(\ref{eq.masks}) 
is the union of the sets of zeros of the functions  
$\bm \bigl( [M^{*}]^{-k}\bxi\bigr)$ over all $k \in \n$
(all zeros are counted with multiplicities). 
Consequently, the set of zeros of the set $\widehat \varphi (\bxi)$ 
is uniquely defined by the set $\bigcup \limits_{k \in \n} M^{-k}\cR$
(all points are counted with multiplicities). Every set $M^{-k}\cR$
lies on the segment $\bigl[-M^{-k}\ba\, , \, M^{-k}\ba\bigr]$, where
$\ba = 2\be_1$ is a point of the set~$\cR$ most distant from the origin, and hence, 
 $M^{-k}\ba$ is a point of the set~$M^{-k}\cR$ most distant from the origin. 
Let $\bar \ba = \bigl[-\ba\, , \, \ba\bigr]$. We see that each set
$\bigcup \limits_{k \in \n} M^{-k}\cR$ is uniquely defined by the set of segments 
$\bigcup \limits_{k \in \n} M^{-k}\bar \ba \, = \, \cT(M, \bar \ba)$. 
Thus, the set of zeros of the function $\widehat \varphi $ 
defines the set of segments  $\cT( M ,\bar \ba)$ in a unique way. 
If we show that this set of segments uniquely, up to a sign, defines the matrix~$M$, 
everything will be proved. If we assume the converse, that there exist integer matrices  $M_1, M_2$ with determinants~$\pm 2$ and segments~$\bar \ba_1$
and $\bar \ba_2$ such that $\cT(M_1, \bar \ba_1)  =  \cT(M_2, \bar \ba_2)$, 
then by Lemma~\ref{l.dyn} we have either $M_1 = \pm M_2$, which is required, 
or one of those matrices has a multiple eigenvalue. The latter is impossible 
by~\cite[Proposition 2.5]{KL00}.

{\hfill $\Box$}
\smallskip

\begin{cor}\label{c.uniq} 
Two polynomials generate affinely similar 2-attractors if and only if they are either 
equal or opposite. 
\end{cor}

In odd dimensions, a polynomial cannot be opposite to itself and every polynomial 
with a free coefficient  $2$ is opposite to a unique polynomial, which has a free coefficient  $-2$, and vice versa. Thus, in odd dimensions, the number of different  (not affinely similar) 2-attractors is equal to the number of integer expanding polynomials with the leading 
coefficient equal to one and the free coefficient equal to~$2$. 

In the even dimensions, opposite polynomials have the same free coefficient. 
Polynomials opposite to themselves are precisely those having zero coefficients with all odd powers. They depend only on $x^2$ and the change $x_1= x^2$ respects the expanding 
property. Therefore, the total number of expanding polynomials opposite to themselves is equal to the number of expanding polynomials of the half degree with the leading coefficient one and the free coefficient~$\pm 2$.  Thus, the number of different (not affinely similar) 
2-attractors in~$\re^{2d}$ is equal to half of the sum of the number for expanding polynomials of degree $2d$ with the leading coefficient one and the free coefficient~$\pm 2$
and of the number for such polynomials of degree~$d$. 
\newpage 

\bigskip 

\begin{center}
\textbf{10. The number of  2-attractors in~$\re^d$}
\end{center}
\bigskip

\begin{figure}[ht]
\centering
\includegraphics[width=0.5\linewidth]{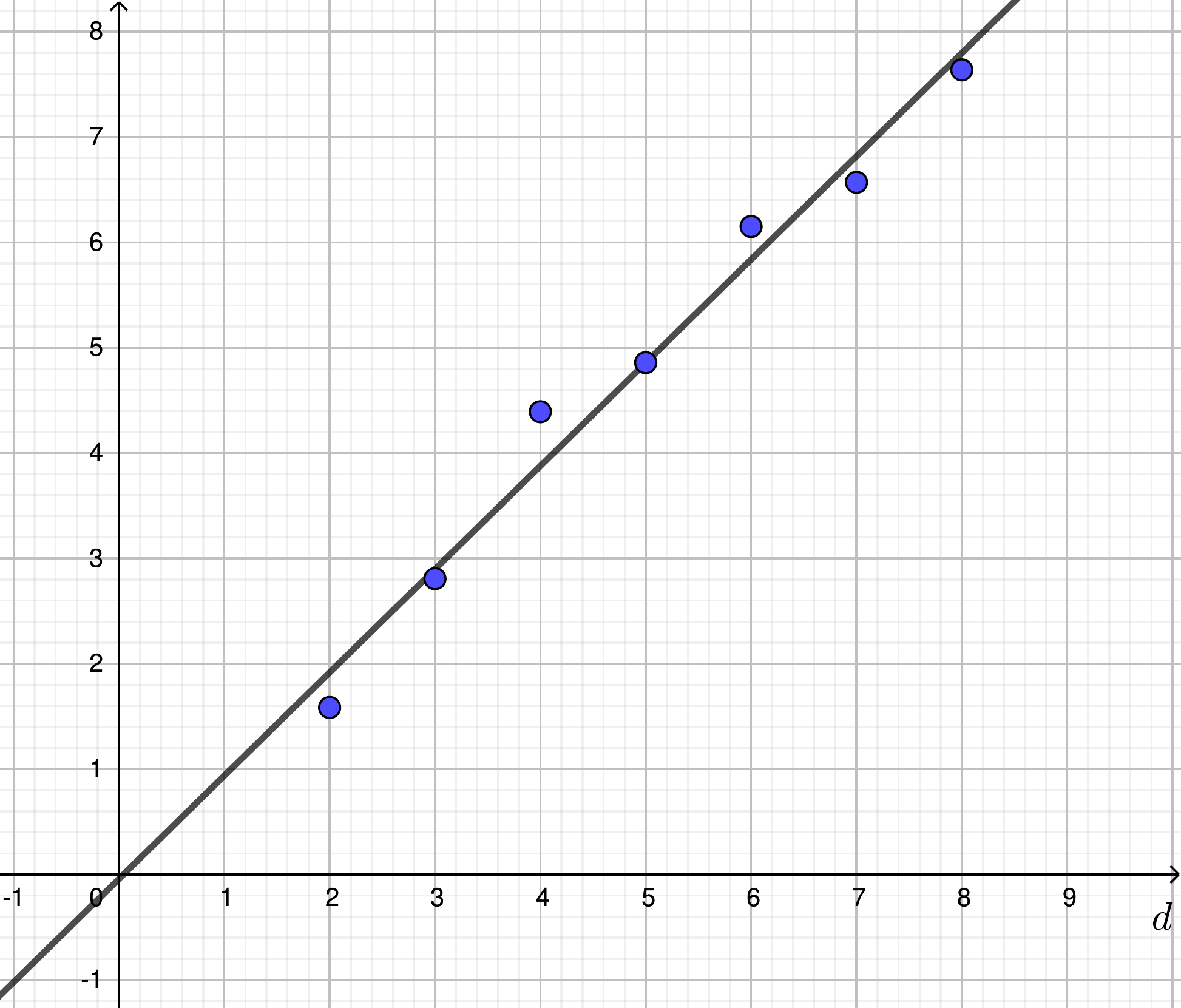}
\caption{The graph of a binary logarithm for the number of 2-attractors in dimension $d$}
\label{graph}
\end{figure}

In every dimension~$d$, there are finitely many 
(not affinely similar) 2-attractors~\cite{Gel}. In $\re^2$ there are exactly three 2-attractors, all of them are isotropic. In  $\re^3$ there are exactly $7\, $ ones, there is only 
one isotropic among them, which is a parallelepiped, see Section~11 for more on the three-dimensional case. Numerical results in Table~\ref{table} show a lower bound for the number of 2-attractors in dimensions $d \le 8$. Most likely, those estimates are quite sharp although we 
are unable to prove this. Figure~\ref{graph} presents a graph of the binary logarithm of this estimate on the number of 2-attractors in dimension~$d$. 
We see that the graph is close to a linear function. 
Therefore, the lower bound for those~$d$ is close to~$2^d$. For bigger dimensions~$d$, 
we know only an approximate answer having obtained lower and upper bounds quite distant from each other.  

\begin{table}[H]
\begin{center}
\begin{tabular}{c|c|c}
d & \text{The number of polynomials} & \text{The number of 2-attractors} \\ \hline
2 & 6 & 3 \\ 
3 & 14 & 7 \\
4 & 36 & 21 \\
5 & 58 & 29 \\
6 & 128 & 71 \\
7 & 190 & 95 \\
8 & 362 & 199 \\
\end{tabular}
\end{center}
\caption{The lower bound for the number of 2-attractors}
\label{table}
\end{table}

By Corollary~\ref{c.uniq} (Section~9), the question on the number of  2-attractors
is reduced to the following problem on integer polynomials: 

\smallskip 

\noindent \textbf{Problem}. {\em How many  integer expanding polynomials of degree $d$
exist with the leading coefficient one and a free term ~$\, \pm \, 2$ ? }
\noindent 
\smallskip 

In~\cite{HSV} it is shown that the number of such polynomials 
grows  at least linearly  in~$d$. On the other hand, this number is finite 
for every $d$. We are going to obtain the following estimates: 
a lower bound is   $\frac{1}{16}\,  d^2 \, + \, O(d)$, and an upper bound is  $\, 2^{\, d+\varepsilon}$. 
Before giving a proof let us remark that the asymptotic distribution of the number of integer expanding polynomials with a large free coefficient has been studied in the literature~\cite{ABPT, AP, KT17}. Polynomials with similar properties appeared also in the work~\cite{Gar}. In \cite{U19}, a simplifying criterion of  expanding property of Schur-Kon for integer 
polynomials was obtained.

We shall start with a proof of the lower bound $\frac{1}{16}\,  d^2 \, + \, O(d)$. 
We call a vector of natural numbers  $(n_1, n_2, n_3)$ \textit{bad} if it is proportional to a vector  $(3x+s, 3y+s, 3z+s)$, where $s \in \{1, 2\}$, and $\, x, y, z \, \in \z$. The corresponding partition of the number  $d=n_1+n_2+n_3$ to a sum of natural terms will be called   \textit{bad} and all other partitions are called  \textit{good}. 
Partitions different  only in the order of the terms will be identified.  By  $b_{+}(d)$ we denote the number of good partitions of the number $d$.

\begin{theorem}\label{t.ser4}
Every polynomial of the form $P(z) = (1 + z^m)(1 + z^q)(1 + z^k) + 1$ 
is an integer expanding polynomial with the leading coefficient~$1$ and the free term~$2$, 
provided $(m, q, k)$ is a good partition $d = m + q  + k$. If  $d$ is not divisible by $3$, 
then the number of such polynomials  $b_{+}(d)$ is equal to~$\frac{d(d-1)}{12} \, + \ \omega_d$, where 
$\omega_d \in \bigl[-\frac23, \frac12 \bigr]$. If $d$ is divisible by~$3$, then it satisfies the equality 
$$
 \frac{d^2}{16}\,  - \, \frac{43d}{36}\, - \frac56 \quad  \le \quad  b_{+}(d)\quad  \le \quad  
 \frac{7d^2}{108}\,  + \, \frac{5d}{12}\, + \frac23\, . 
$$
\end{theorem}

Note that the principal terms of the upper bound and of the lower bound for  $b_{+}(d)$ are 
quite  close: 
$\frac{1}{16} = 0.0625, \frac{7}{108} \approx 0.0648$. The proof of Theorem~\ref{t.ser4}
is postponed to Section~12. 
\smallskip 

To derive the upper bound $2^{\, d+\varepsilon}$ we apply a result of A.Dubickas and S.Konyagin 
on the number of polynomials whose Mahler measure does not exceed a number~$T$.

The {\em Mahler measure} of an algebraic polynomial~$p(t)\, = \, a_dt^d + \cdots + a_1t + a_0$ is the number 
$\mu(p) \, = \, |a_d|\, \prod_{i=1}^d \max\, \{1, |\lambda_i|\}$, where $\lambda_1, \ldots , \lambda_d$ are roots of~$p$ counting multiplicity. The following theorem from~\cite{DK}
will be cited in a slightly simplified form. We denote by $\theta = 1.32471...$ the root of the polynomial $x^3 - x - 1$.  
\smallskip 

\noindent \textbf{Theorem B}~\cite{DK}. If $T \ge \theta$, then the number of all integer polynomials with the leading coefficient~$1$ 
and with the Mahler measure at most~$T$ does not exceed  $T^{\, d\, \bigl(1+\frac{16 \ln \ln d}{\ln d}\bigr)}$. 
\noindent 
\smallskip 
 
Now we are ready to formulate the main theorem of this section. 
\begin{theorem}\label{th.estimates}
The total number  $N(d)$ of not affinely similar 2-attractors in dimension~$d$
satisfies the following estimates: $\quad \frac{d^2}{16}\,  - \, \frac{43d}{36}\, - \frac56 \, \le \, N(d) \, \le \, 2^{\, d\, \bigl(1+\frac{16 \ln \ln d}{\ln d}\bigr)}$  
\end{theorem}
{\tt Proof}.  
 The lower bound follows from Theorem~\ref{t.ser4}, which estimates the number of not opposite 
 each other polynomials of the form $P(z) = (1 + z^m)(1 + z^q)(1 + z^k) + 1$. By 
 Theorem~\ref{th.uniq}, all of them define not affinely similar 2-attractors.

 To establish the upper bound we observe that the Mahler measure of an expanding polynomial 
 with coefficients  $a_d = 1, a_0 = \pm 2$ is equal to  $\prod_{i=1}^d |\lambda_i| \, = \, |a_0|\, = \, 2$. 
Hence, $N(d)$ does not exceed the number of integer polynomials satisfying $\mu(p) \le 2$. 
By applying Theorem~B for $T=2$ we complete the proof.  

{\hfill $\Box$}
\smallskip

\smallskip

A large difference between the lower and upper bounds in Theorem~\ref{th.estimates}
rises the question on the asymptotic behaviour of the value~$N(d)$ as $d\to \infty$. 
At the first sight,~$N(d)$ has to be significantly smaller than the number of polynomials with Mahler's measure~$\mu(p) \le 2$. First, because in $N(d)$, we count only 
polynomials with the leading coefficient~$1$. Second, because roots smaller than one by modulus are prohibited. Nevertheless, the numerical results from Table~\ref{table} 
rather show the opposite: for $d\le 8$ the ratio  $N(d)/2^d$ approaches to one when~$d$
grows. Observe also that lower bounds on the number of polynomials with a prescribed Mahler measure are often significantly smaller that the upper bound since examples of such polynomials 
``breed poorly'', see~\cite{Dub}. 

\bigskip

\begin{center}
\textbf{11. 2-attractors in~$\re^3$}
\end{center}
\bigskip 

In the three-dimensional space there exist precisely $7$ types of non affinely similar 2-attractors. They were first constructed in~\cite{BG}. 
Those 2-attractors have characteristic polynomials different and not opposite to each other. 
 Hence, by Corollary~\ref{c.uniq}, they are not affinely similar. 
 The work~\cite{B10} studies their topological properties: 
 homeomorphism to a ball, the structure of neighbours, etc. 
 In particular, the combinatorial structure of tilings of~$\re^3$ by integer shifts of those attractors was found in that work. In this section we compute the H\"older exponents of three-dimensional attractors and will see that all of them are different. 
 This explains not only their affine non-equivalence but also  different 
 approximation properties of the corresponding Haar wavelets.
We shall see that even drawing pictures of attractors significantly depends on their 
smoothness.  

Recall that every polynomial $\ p(z) \, = \, z^3 + a_{2}z^2 + a_1z + a_0$ 
is associated to its companion matrix 
\begin{equation*}
M \quad = \quad \left(
\begin{array}{ccc}
0 &  0 & -a_0\\
1 &  0 & -a_{1} \\
0 &  1 & -a_{2} \\
\end{array}
\right) 
\end{equation*}
with the characteristic polynomial~$p$.
For each of the seven polynomials written in Table~\ref{sm_table} we build an attractor 
with the companion matrix and with the digits  $(0, 0, 0)^T$ and $(1, 0, 0)^T$. 
The H\"older exponent  $\alpha (G)$ of a set $G$ is the H\"older exponent  in $L_2$ of its characteristic function $\chi_G$ (see Definition~(\ref{eq.Holder})). 

The H\"older exponent is evaluated using the method from~\cite{CP}. 
First we consider auxiliary matrices $T_0$, $T_1$, their dimensions depend on the attractor
 (from $4$ to $23$ in our seven cases). Then we compute the $L_2$-spectral radius~\cite{P97} of those matrices on a special subspace, which gives the H\"older exponent. 

\begin{table}[H]
\begin{center}
\begin{tabular}{c|c|c}
\text{The polynomial} & \text{$L_2$-spectral radius} & \text{H\"older exponent in $L_2$} \\ \hline
$z^3 + 2z^2 + 2z + 2$ & 0.97082 & 0.06822 \\ 
$z^3 + z^2 + z + 2$ & 0.93238 & 0.23148 \\
$z^3 + 2$ & 0.8909 & 0.5 \\
$z^3 - z^2 - z + 2$ & 0.94278 & 0.23282 \\
$z^3  - z + 2$ & 0.95197 & 0.1173 \\
$z^3 - 2z + 2$ & 0.98548 & 0.02563 \\
$z^3 + z^2 + 2$ & 0.97542 & 0.04713 \\
\end{tabular}
\end{center}
\caption{The regularity of three-dimensional 2-attractors}
\label{sm_table}
\end{table}

\begin{figure}[ht!]
\begin{minipage}[h]{0.3\linewidth}
\center{\includegraphics[width=1\linewidth]{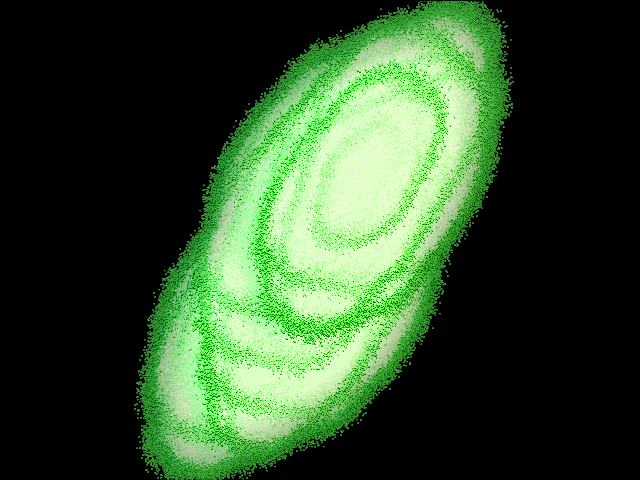}} 1 a) \\
\end{minipage}
\hfill
\begin{minipage}[h]{0.3\linewidth}
\center{\includegraphics[width=1\linewidth]{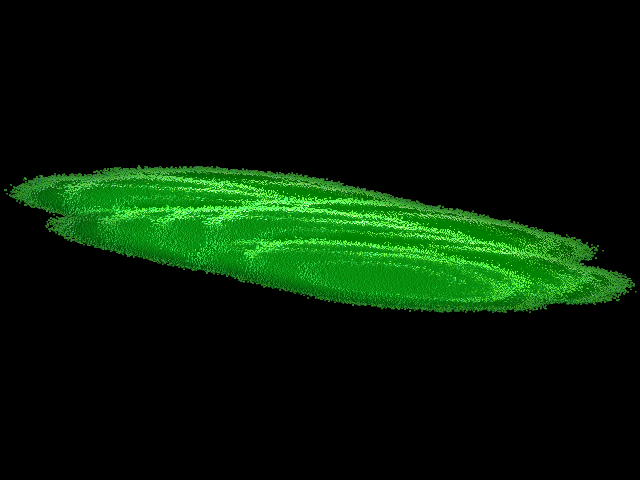}} 1 b)\\
\end{minipage}
\hfill
\begin{minipage}[h]{0.3\linewidth}
\center{\includegraphics[width=1\linewidth]{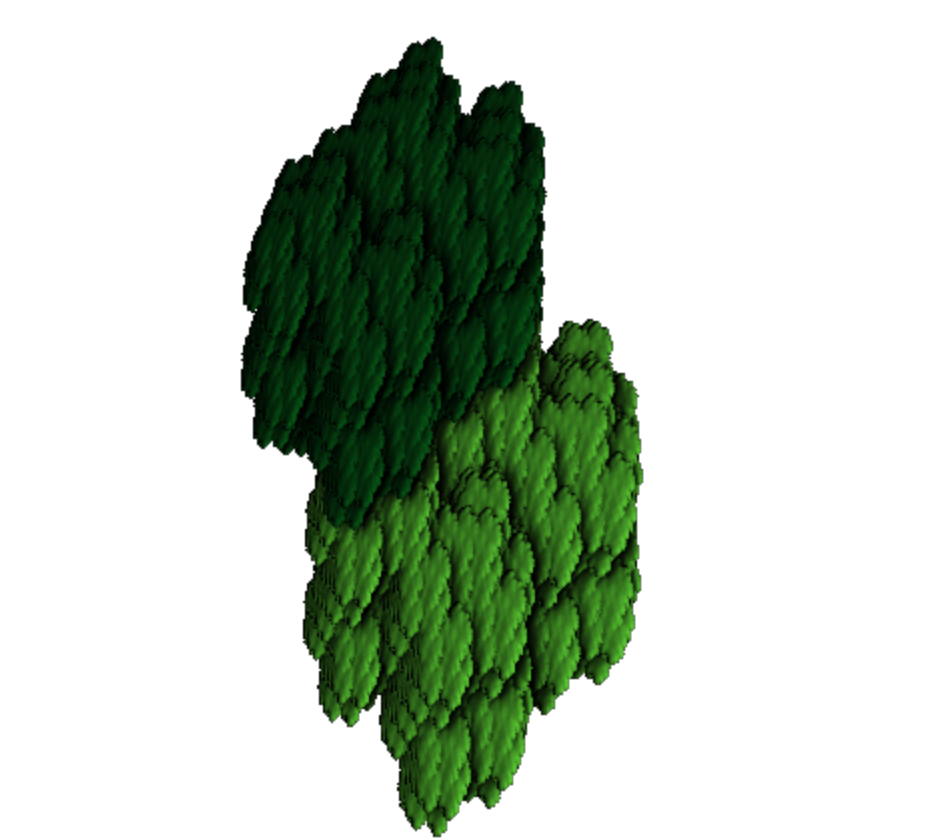}} 2 a)\\
\end{minipage}

\vfill

\begin{minipage}[h]{0.3\linewidth}
\center{\includegraphics[width=1\linewidth]{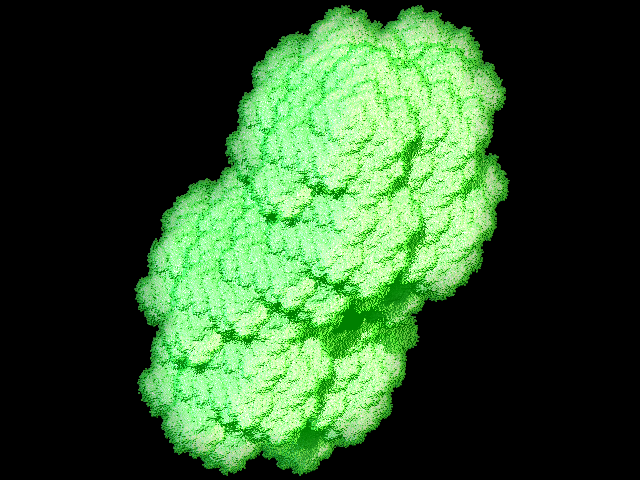}} 2 b) \\
\end{minipage}
\hfill
\begin{minipage}[h]{0.3\linewidth}
\center{\includegraphics[width=1\linewidth]{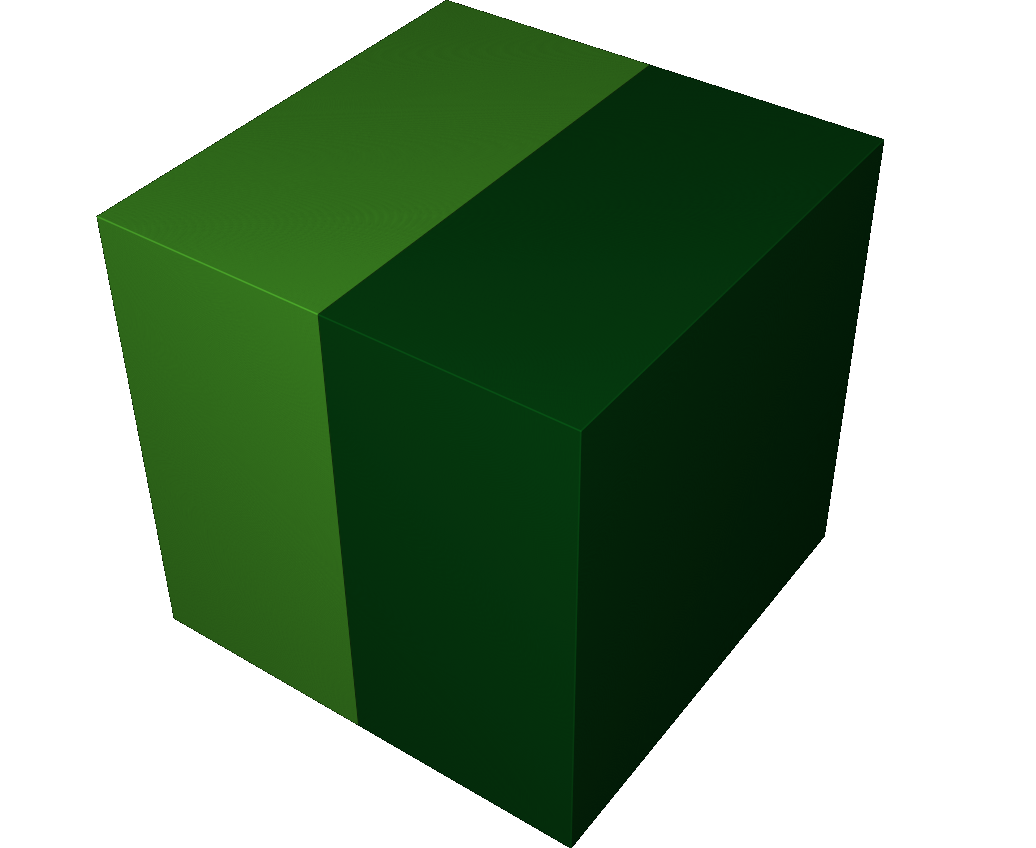}} 3 \\
\end{minipage}
\hfill
\begin{minipage}[h]{0.3\linewidth}
\center{\includegraphics[width=1\linewidth]{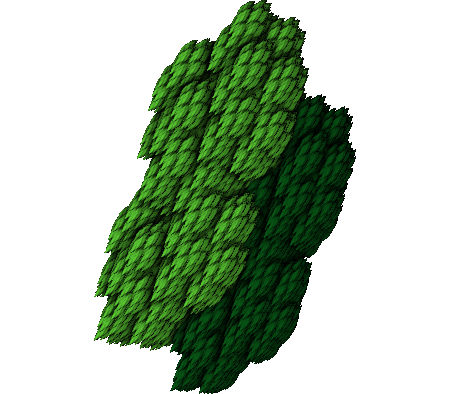}} 4 a) \\
\end{minipage}
\caption{The seven types of three-dimensional 2-attractors, part 1}
\label{pic7_1}
\end{figure}

\begin{figure}[ht!]
\begin{minipage}[h]{0.3\linewidth}
\center{\includegraphics[width=1\linewidth]{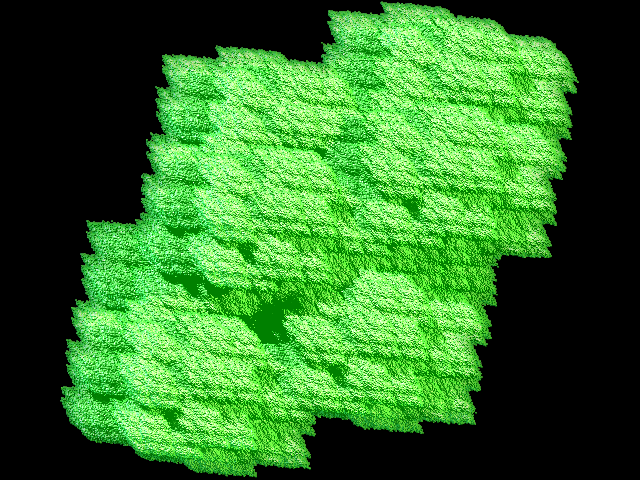}} 4 b) \\
\end{minipage}
\hfill
\begin{minipage}[h]{0.3\linewidth}
\center{\includegraphics[width=1\linewidth]{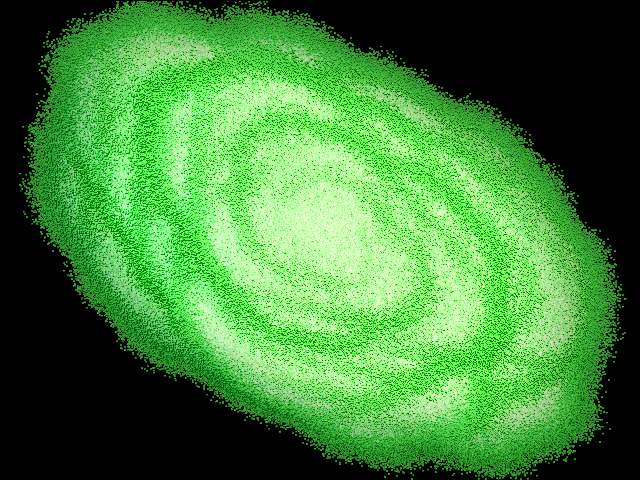}} 5 a)\\
\end{minipage}
\hfill
\begin{minipage}[h]{0.3\linewidth}
\center{\includegraphics[width=1\linewidth]{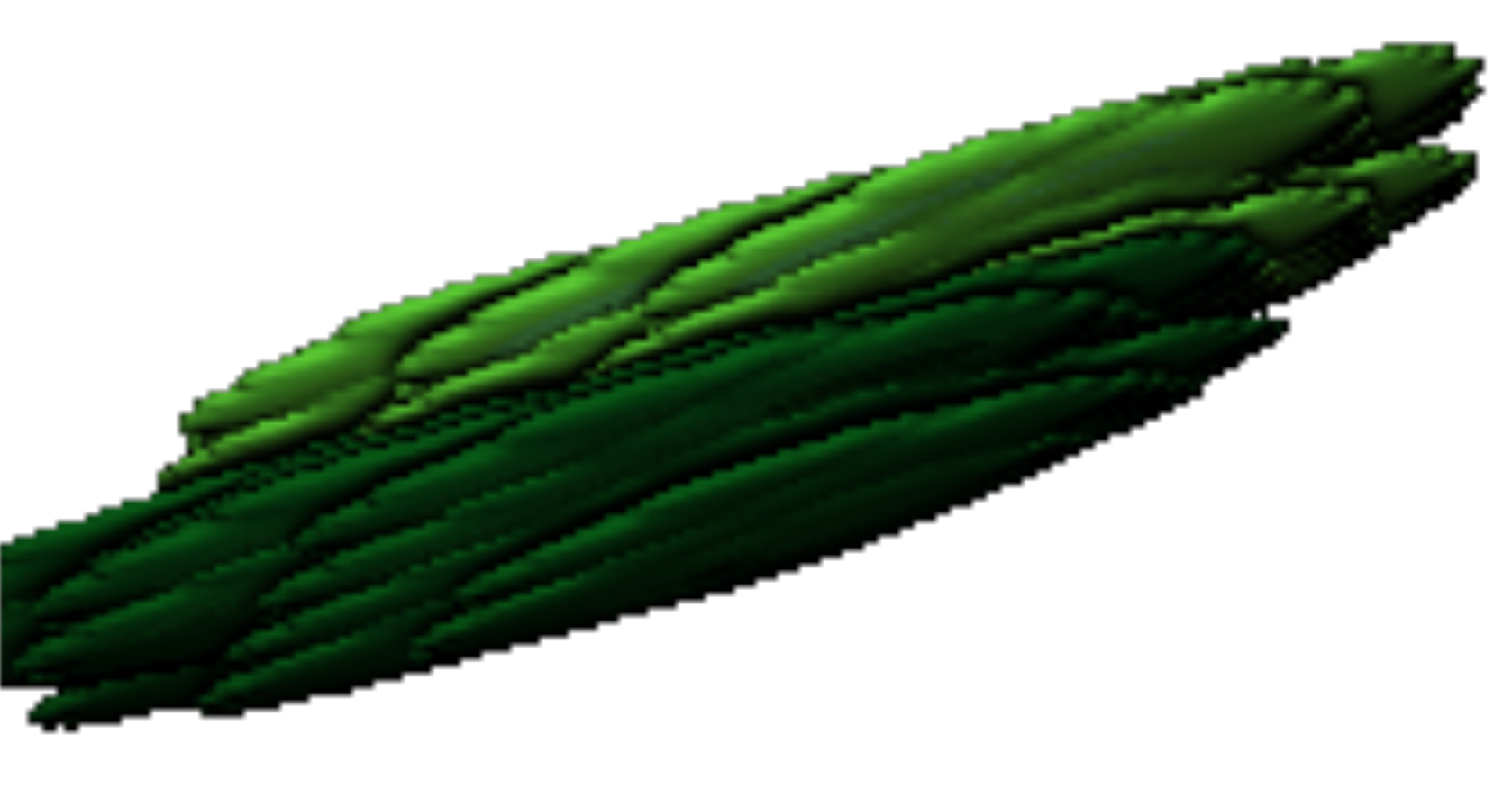}} 5 b)\\
\end{minipage}

\vfill

\begin{minipage}[h]{0.3\linewidth}
\center{\includegraphics[width=1\linewidth]{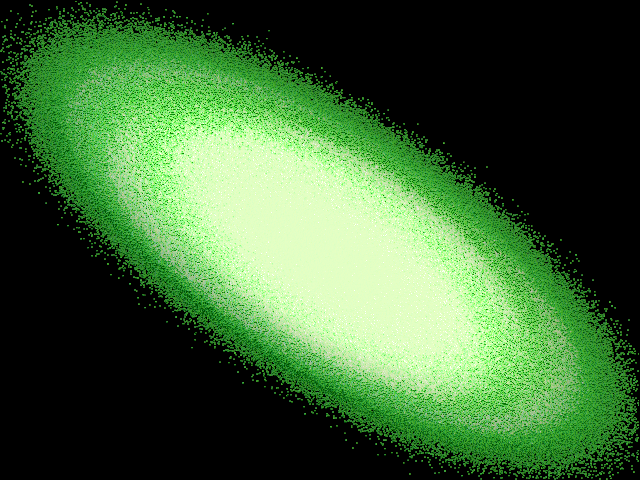}} 6 a) \\
\end{minipage}
\hfill
\begin{minipage}[h]{0.3\linewidth}
\center{\includegraphics[width=1\linewidth]{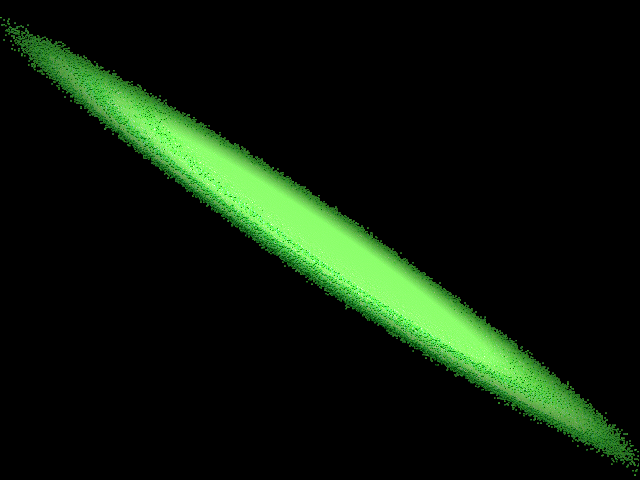}} 6 b) \\
\end{minipage}
\hfill
\begin{minipage}[h]{0.3\linewidth}
\center{\includegraphics[width=1\linewidth]{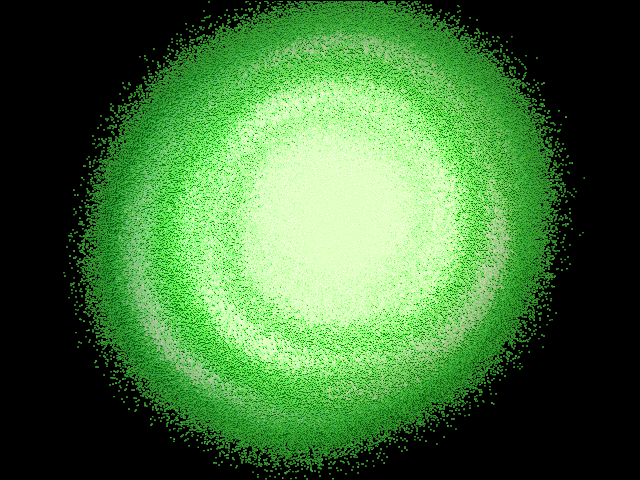}} 7 \\
\end{minipage}
\caption{The seven types of three-dimensional 2-attractors, part 2}
\label{pic7_2}
\end{figure}

The two-dimensional illustrations of 2-attractors  (with the dragon, rectangle, and the bear) are drown with the program~\cite{KM}. Three-dimensional illustrations are made with two programs, Chaoscope \cite{Desp} and  IFStile \cite{Mekh}. Their algorithms are different.  Chaoscope uses the probabilistic approach, the running time slightly depends on the type of attractor. 
The program IFStile producing more detailed pictures uses,  most likely, 
the iteration principle. Its running time significantly depends on 
the H\"older regularity of the attractor~(see Remark~\ref{r.quality-draw} below). 
For example, the construction of the third type attractor  (parallelepiped), which has 
the maximal regularity $\alpha = 0.5$, takes less than a minute, 
the type 2 ($\alpha \approx 0.2315$) and type 4 ($\alpha\approx0.2328$) 
with the same parameters require several hours, the type 5 ($\alpha \approx 0.1173$) 
needs already several days, even with a worse image quality. 
Those pictures also show the partition of each attractor into two 
affinely similar parts  (dark green and light green). The images are in Figures~\ref{pic7_1} and~\ref{pic7_2}, some attractors are shown from different angles. 

\begin{remark}\label{r.quality-draw}
{\em The quality of the attractors from the Figures~\ref{pic7_1}, \ref{pic7_2} significantly depends on their regularity. Let us try to explain this phenomenon by the example of the iteration method, which, 
probably, more or less, underlies all existing algorithms.
For the linear operator $T$ acting in  $L_2(\re^d)$ by the formula 
$$
[Tf](\bx) \ = \ \sum_{\bk \in D}\, f(M\bx - \bk)\, ,
$$
the characteristic function $\varphi = \chi_{\, G}$ is a fixed point. 
If the attractor $G$ is a tile, then for the initial function 
$f_0 = \chi_{[0,1]^d}$ (the characteristic function of the unit cube),  
the iteration algorithm $f_{k+1} = Tf_{k}, k \ge 0$,  converges to~$\varphi$ in $L_2$. 
Moreover, $\|f_k - \varphi\|_{L_2} = O(2^{- k \alpha})$, 
where $\alpha$ is the H\"older exponent of $\varphi$ (this fact is known for every 
refinable function~\cite{CDM, CHM, KPS}). Therefore,   
$k \, \approx \, \frac{\log_2 \frac{1}{\varepsilon}}{\alpha}$. Say, for 
approximating  $\varphi$ with precision $0.03$ one needs about $k = \frac{5}{\alpha}$ 
iterations. For the attractor of type~3, one needs 10 iterations 
 (less than a minute of the computer time). 
For the types 2 and 4 we need already  22 iterations, which requires several hours, 
because the complexity of one iteration increases. The type 5 needs  
 $45$ iterations, which is hardly realisable within a reasonable time without an essential 
 lost of the image quality.  The types  6 and 7 require 195 (!) and 106 iterations respectively.

}
\end{remark}

\bigskip 

\begin{center}
\textbf{12. Series of integer expanding polynomials and 2-attractors}
\end{center}
\bigskip 

As we know (Corollary~\ref{c.uniq}),  a 2-attractor is uniquely defined by an integer
expanding polynomial~$p(z)$ with the leading coefficient~$1$ and a free term  $\pm 2$. 
To construct an attractor, one needs to take an integer matrix~$M$ with characteristic 
polynomial~$p$ (for example, its companion matrix~(\ref{eq.accomp}) or any other) 
and an arbitrary proper pair of digits~$D$. The obtained 2-attractor will be 
defined in a unique way, up to an affine similarity, independently of~$M$ and $D$.   
Hence, constructing 2-attractors is equivalent to finding the corresponding 
expanding polynomial. In this section we find several series of such polynomials. 
The number of polynomials of degree~$d$ in each  series grows linearly in~$d$, 
only  series~\ref{ser4} contains a quadratic in $d$ number of polynomials.  
This is that series providing the lower bound on the number of 2-attractors in Theorem~\ref{th.estimates}. 
 We first present all series and then give proofs. Series~\ref{ser1} and \ref{ser2} 
 are known in the literature, series \ref{ser3} -- \ref{ser7}, 
to the best of our knowledge, are new. 
\smallskip 

\textbf{Series of integer expanding polynomials with the leading coefficient one and the free term  $\pm 2$}

\smallskip 
\begin{ser} \label{ser1}
a) Polynomials of the form $\frac{z^m + z^q - 2}{z^{(m, q)} - 1}$ with $m > q \ge 1$, 
where $(m, q)$ is the greatest common divisor (g.c.d) of $m$ and $q$.  

b) $\frac{z^m + z^q + 2}{z^{(m, q)} + 1}$ whenever $m_1 = \frac{m}{(m, q)}$ and $m_2 = \frac{q}{(m, q)}$ are odd, $m > q \ge 1$. 
\end{ser}

\begin{ser} \label{ser2}
Polynomials of these series have only three terms:   

a) $z^m - z^q + 2$ with  $m > q \ge 1$ and  odd $q_1 = \frac{q}{(m, q)}$. 

b) The opposite series $z^q - z^m - 2$ with $q > m \ge 1$ and odd $q_1 = \frac{q}{(m, q)}$. 

c) $z^m + z^q + 2$ with $m > q \ge 1$  whenever the numbers $m_1 = \frac{m}{(m, q)}$ and $q_1 = \frac{q}{(m, q)}$ have different evenness. 
\end{ser}

\begin{ser} \label{ser3}
a)  $(1 + z^m)(1 + z^q) + 1$ for all natural $m$, $q$. 

b) $(z^m - 1)(z^q - 1) + 1$ for all natural $m$, $q$. 
\end{ser}
\smallskip 

Let us remember that a vector with natural components  $(n_1, n_2, n_3)$ 
is called bad if it is proportional to a vector of the form $(3x+s, 3y+s, 3z+s)$ with some  $s \in \{1, 2\}$ and  $\, x, y, z \, \in \z$. All other vectors are called good.

\begin{ser}\label{ser4}
a)  $(1 + z^m)(1 + z^q)(1 + z^k) + 1$ for an arbitrary good vector  $(m, q, k)$. 

b) $(1 - z^m)(1 - z^q)(1 + z^k) + 1$  apart from the case when
$m$ and  $q$ are equal modulo $3$ and $k \equiv -m \pmod{3}$.
\end{ser}

\begin{ser} \label{ser5}
$(1 + z^m + z^{2m})(1 + z^q + z^{2q}) + 1$ whenever $m \not \equiv q \pmod{4}$. 
\end{ser}

\begin{ser} \label{ser6}
$(1 + z^m)(1 \pm z^r + z^{2r}) + 1$ with arbitrary natural  $m$, $r$. 
\end{ser}

\begin{ser} \label{ser7}
$(1 + z^a)(1 + z^b)(1 + z^{(a + b)}) (1 + z^{2(a + b)}) (1 + z^{4(a + b)}) (1 + z^{8 (a + b))}) \ldots (1 + z^{2^k (a + b)}) + 1$, where $a$, $b$, $k$ are natural numbers. 
\end{ser}

\begin{remark}{\em 
It is easy to see that polynomials of series  \ref{ser1} -- \ref{ser7} are different 
(apart from the case $k = q$ in \ref{ser4} b), which is reduced to \ref{ser3} b) and from 
special cases of some series that are covered by \ref{ser1}). All polynomials of series \ref{ser2} a) b) c), \ref{ser3} a) b), \ref{ser4} a), \ref{ser5}, \ref{ser6}, \ref{ser7} take different values at the point $z = 1$: respectively $2$, $-2$, $4$, $5$, $1$, $9$, $10$, $3$ or $7$, $2^{k + 3} + 1$ (where $k \le 1$). }
\end{remark}
\smallskip 

Now we prove that these polynomials are expanding, i.e, all their roots are out of the closed 
unit disc in the complex plane. 
\smallskip  

{\tt Proof for series~\ref{ser1}}.   
a) Let $m = (m, q)m_1$, $q = (m, q)q_1$. The polynomial can be rewritten as a sum 
of two geometric progressions 
$$\frac{z^m - 1}{z^{(m, q)} - 1} + \frac{z^q - 1}{z^{(m, q)} - 1} = 1 + {z^{(m, q)}} + {(z^{(m,q)})}^2 + \ldots + {(z^{(m,q)})}^{m_1 - 1} + 1 + {z^{(m, q)}} + \ldots + {(z^{(m,q)})}^{q_1 - 1}.$$ 
Therefore, the free term is equal to $2$ and the leading coefficient is one. 
If the denominator vanishes: $z^{(m, q)} - 1 = 0$, then the value of the polynomial 
is equal to $m_1 + q_1$ and hence, such  $z$ are not roots of the polynomial. 
In the sequel we assume that $z^{(m, q)} \ne 1$. Then $z^m + z^q - 2 = 0$,  
hence, the roots cannot be smaller than one by modulus. 
If the modulus of a root is equal to one, then the triangle inequality yields  $z^m = z^q = 1$.  Denote by  $\gamma$ the argument of $z$ chosen on the interval $[0, 2\pi)$. 
Since $z^m = 1$, we have $\gamma = \frac{2\pi m_2}{m}$, $m_2 = 0, \ldots, m - 1$. 
Moreover, $\gamma = \frac{2\pi q_2}{q}$, $q_2 = 0, \ldots, q - 1$, so  $\frac{2\pi m_2}{m} = \frac{2\pi q_2}{q}$, which implies  $m_2 q = q_2 m$ and then  $m_2 q_1 = q_2 m_1$. 
Since $(m_1, q_1) = 1$, we see that $m_2$ is a multiple of $m_1$. However, in this case, the polar angle of $z^{(m, q)}$ is equal to $\frac{2\pi m_2 (m, q)}{m_1 (m, q)}$ ,  
which is a multiple of $2\pi$. Therefore, the denominator $z^{(m, q)} - 1$ vanishes, which contradicts to the assumption. So, this polynomial does not have roots such that  $|z| \le 1$. 
\smallskip 

b) This is a similar series and can be obtained by changes the signs. In this case the polynomial has the form
$$\frac{z^m + 1}{z^{(m, q)} + 1} + \frac{z^q + 1}{z^{(m, q)} + 1} = 1 - {z^{(m, q)}} + {(z^{(m,q)})}^2 + \ldots + {(z^{(m,q)})}^{m_1 - 1} + 1 - {z^{(m, q)}} + \ldots + {(z^{(m,q)})}^{q_1 - 1}.$$
It can be assumed again  $z^{(m, q)} + 1 \ne 0$, otherwise  $z$ is not a root. 
Again there are no roots inside the unit circle, and on the unit circle we have  $z^m = -1$, $z^q= -1$. 
In this case the argument of the number $z$ equal to $\gamma \in [0, 2\pi)$ satisfies the equalities  $\gamma = \frac{\pi}{m} + \frac{2\pi m_2}{m} = \frac{\pi}{q} + \frac{2\pi q_2}{q}$. After the multiplication by  $\frac{mq}{\pi (m, q)}$ we obtain $q_1 + 2m_2q_1 = m_1 + 2q_2m_1$. Since $(q_1, m_1) = 1$, it follows that  $(1 + 2m_2)$ is divisible by $m_1$. 
The polar angle of the number $z^{(m, q)}$ is   $\frac{\pi(m, q)}{m} + \frac{2\pi m_2(m, q)}{m} = \pi \frac{1 + 2m_2}{m_1}$. This point corresponds to the number $-1$ because $\frac{1 + 2m_2}{m_1}$ is an odd number. The denominator vanishes,  which contradicts to the assumption.

{\hfill $\Box$}
\smallskip

{\tt Proof for  series \ref{ser2}}.   
In all the three cases the roots cannot be inside the unit circle since 
in this case $|z^m \pm z^q| < 2$. 

If $|z| = 1$, then in  cases a) and b) we have $z^m = -1$, $z^q = 1$, and in  case c) we have $z^m = -1$, $z^q = -1$. 

In cases a) and b) we get  $\frac{\pi}{m} + \frac{2\pi m_2}{m} = \frac{2\pi q_2}{q}$, where $m_2 = 0, \ldots, m - 1$, $q_2 = 0, \ldots, q - 1$. After multiplying by  $\frac{mq}{\pi (m, q)}$ we obtain  $q_1 (1 + 2m_2) = 2q_2 m_1$. This is impossible since  $q_1$ and $1 + 2m_2$ are both odd. 

In case c) we have   $\frac{\pi}{m} + \frac{2\pi m_2}{m} = \frac{\pi}{q} + \frac{2\pi q_2}{q}$, therefore,  $q_1+ 2m_2 q_1 = m_1 + 2m_1q_2$. Hence, $q_1 - m_1 = 2 (m_1q_2 - m_2q_1)$, which is impossible provided $q_1$ and $m_1$ are of different evenness.

{\hfill $\Box$}
\smallskip

{\tt Proof for series \ref{ser3}}.   
a) We assume the converse: our polynomial $P(z)$ has a root $|z| \le 1$.
 
Denote by  $\overline{B}(1, 1)$ the disc in the complex plane centered at  $1$ and of radius $1$. 
Clearly, $|z^m| \le 1$ and $|z^q| \le 1$, therefore, both points $z^m + 1$ and $z^q + 1$
lie in the disc  $\overline{B}(1, 1)$. Observe that $z^m + 1 \ne 0$ and  $z^q + 1 \ne 0$, 
otherwise $P(z) = 1 \ne 0$. Every nonzero complex number in the disc $\overline{B}(1, 1)$
has the argument strictly less than $\frac{\pi}{2}$ by modulus. Hence, the product of 
two arbitrary nonzero numbers from that disc has the argument less than $\pi$ by modulus. Therefore, the product of the numbers $z^m + 1$ and $z^q + 1$ 
cannot have the argument $\pi$ and be equal to $-1$, which is a contradiction. 

b) The proof for b) is literally the same.

{\hfill $\Box$}
\smallskip

{\tt Proof for series \ref{ser4}}.   
a) Assume our polynomial $P(z)$ has a root $|z| \le 1$. Consider points $z_1 = z^m + 1$, $z_2 = z^q + 1$, $z_3 = z^k + 1$, they belong to the disc $\overline{B}(1, 1)$. 
Denote their arguments by $\alpha_1, \alpha_2, \alpha_3$ respectively. Note that $z_i \ne 0, i = 1, 2, 3$, otherwise $P(z) = 1 \ne 0$. Each interior point of the disc $\overline{B}(1, 1)$ has an argument smaller than  $\frac{\pi}{2}$. Hence, $|\alpha_i| < \frac{\pi}{2}, i= 1, 2, 3$. Since  $z_1z_2z_3 = -1$, it follows that $\alpha_1 + \alpha_2 + \alpha_3 = \pm \pi$. 
This implies that all the numbers $\alpha_i$ have the same sign. 
It can be assumed that they are positive, otherwise we replace $z$ by the conjugate root. 
Thus, the angles  $\alpha_i$ are angles of an acute triangle. 
For an acute triangle, the following inequality holds:  $\cos(\alpha_1)\cos(\alpha_2)\cos(\alpha_3) \le \frac{1}{8}$, it becomes an equality only for an equilateral triangle. 
For convenience of the reader, we prove this assertion.   The function $f(x) = -\ln(\cos(x))$  is strictly convex on the interval $x \in (0, \frac{\pi}{2})$ because $f''(x) = \frac{1}{\cos^2(x)} > 0$. We need to find a maximum of the function $\cos(\alpha_1)\cos(\alpha_2)\cos(\alpha_3)$, or, which is equivalent, to minimize the function  $f(\alpha_1, \alpha_2, \alpha_3) = -\ln(\cos(\alpha_1)) - \ln( \cos(\alpha_2)) - \ln(\cos(\alpha_3))$, 
which is coercive on a compact set  $[0, \frac{\pi}{2}]^3$ (convex with a possible value $+\infty$). Such a function possesses at most one point of minimum, which is attained on the intersection with the hyperplane  $\alpha_1 + \alpha_2 + \alpha_3 = \pi$ when all the arguments equal to each other. 

Thus, $|\cos(\alpha_1)\cos(\alpha_2)\cos(\alpha_3)| \le \frac{1}{8}$. Since $z_i \in \overline{B}(1, 1)$, we have $|z_i| \le 2 \cos(\alpha_i)$. On the other hand $|z_1 z_2 z_3| = 1$, this implies  $|\cos(\alpha_1)\cos(\alpha_2)\cos(\alpha_3)| \ge \frac{1}{8}$. We come to the only possible case when the triangle is equilateral and  $|z_1| = |z_2| = |z_3| = 1$. 
Then  $z^m + 1$, $z^q + 1$, and $z^k + 1$ have the same argument $\frac{\pi}{3}$ and hence,  $z^m$, $z^q$, and $z^k$ have the same argument $\frac{2\pi}{3}$.  We denote the 
argument of the number  $z$ by $\gamma \in [0, 2\pi)$; then $\gamma = \frac{2\pi}{3m} + \frac{2\pi m_2}{m}= \frac{2\pi}{3q} + \frac{2\pi q_2}{q}$, $m_2 = 1, \ldots, m - 1$, $q_2 = 1, \ldots, q - 1$. By multiplying by $\frac{3mq}{2\pi}$, we obtain $q + 3m_2q = m + 3q_2m$. 
Hence, the numbers  $q$ and $m$ have the same remainders upon division by $3$. 
The number $k$ must have the same reminder. If these reminders are not zero,  then we have already obtained a bad vector $(q, m, k)$. If all the numbers $(q, m, k)$ are divisible by $3$, then we divide them by the maximal common power of $3$. We obtain numbers with the same 
nonzero remainder upon division by $3$, i.e., in this case $(q, m, k)$ is also a bad vector. Hence, for good vectors, the polynomial is expanding.

b) The proof remains the same with the change of notation $z_1 = 1 - z^m$, $z_2 = 1 - z^q$, 
until we come to the equalities $|z_1| = |z_2| = |z_3| = 1$ and $\alpha_1 = \alpha_2 = \alpha_3 = \frac{\pi}{3}$ (we again without loss of generality assume that the angles are positive). 
However, this time the argument of $z^m$ and of $z^q$ is  $\frac{4\pi}{3}$, and 
the argument of $z^k$ is again  $\frac{2\pi}{3}$. Denote the argument of  $z$ by $\gamma \in [0, 2\pi)$, then $\gamma = \frac{4\pi}{3m} + \frac{2\pi m_2}{m}= \frac{2\pi}{3k} + \frac{2\pi k_2}{k}$, $m_2 = 1, \ldots, m - 1$, $k_2 = 1, \ldots, k - 1$. Multiply by $\frac{3mk}{2\pi}$
and obtain $2k + 3m_2k = m + 3q_2m$. Hence, $m \equiv -k \pmod{3}$. Analogously,  $q \equiv -k \pmod{3}$. If these conditions fail, then the polynomial is expanding.

{\hfill $\Box$}
\smallskip

{\tt Proof for series \ref{ser5}}.   
Denote our polynomial by $P(z) = (1 + z^m + z^{2m})(1 + z^q + z^{2q}) + 1$. 
We prove that if $(1 + z^m + z^{2m})(1 + z^q + z^{2q}) = - c$, where $c$ is a real positive number and  $|z| \le 1$, then $c < 1$. This is sufficient, since in this case 
$P(z)$ does not have roots such that  $|z| \le 1$. Since the polynomial is an analytic function, 
it follows by the maximum modulus principle that it suffices to prove the assertion in the case 
  $|z| = 1$, i.e., for $z = e^{i \alpha}$. 

\begin{figure}[ht]
\centering
\includegraphics[width = 0.6\textwidth]{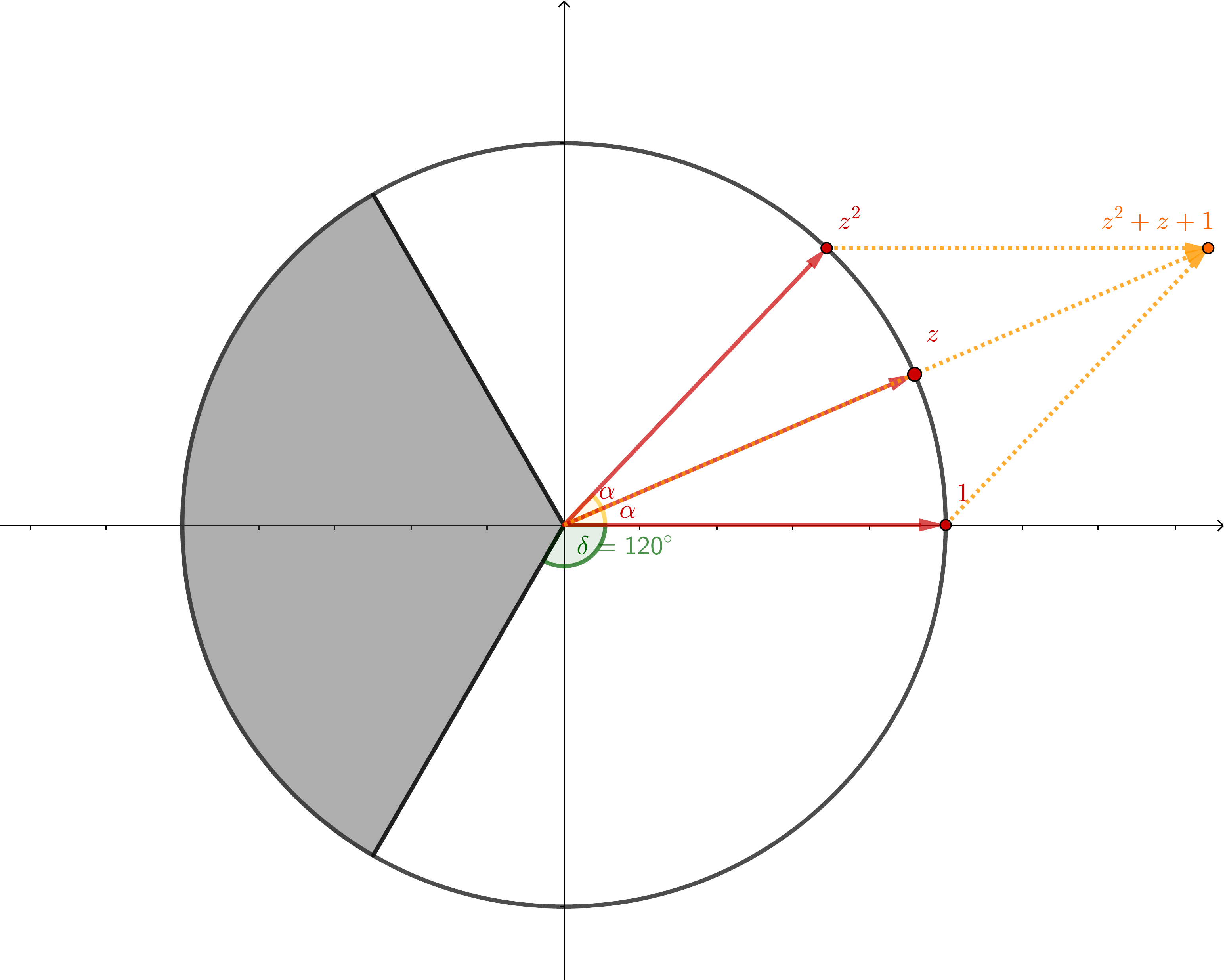}
\caption{The sum of the form $1 + z + z^2$}
\label{zones}
\end{figure}

Let us investigate the properties of the polynomial $p(z) = 1 + z + z^2$ for $|z|= 1$. We spot a ``white zone'' on the complex plane which consists of points on the unit circle with angles on 
  $[-\frac{2\pi}{3}, \frac{2\pi}{3}]$. Its complement on the circle will be called a ``grey zone'' (see Figure~\ref{zones}). Observe that the summation of the vectors  $z^2$ and $1$ 
  gives a parallelogram with diagonal parallel to $z$. If the point  $z$ is in the white zone, then  $p(z) = \lambda z$ for $\lambda \ge 0$, if $z$ is in the grey zone, then $p(z) = \lambda z$ for $\lambda < 0$. In both cases we obtain that  $p(z)$ is in the white zone. 
Furthermore, we have  $|p(z)| = |1 + z + z^2| =  |1 + 2\cos \alpha|$. 

We consider the question, for which positive real $c$ 
 we can have $p(u)p(w) = -c$, where $u = z^m$, $w = z^q$, $|u| = |w| = |z| = 1$. If $u$ 
is in the grey zone, then the argument  $p(u)$ lies in $(-\frac{\pi}{3}, \frac{\pi}{3})$. The argument $p(w)$ is in $[-\frac{2\pi}{3}, \frac{2\pi}{3}]$, since $p(w)$ is in the white zone. Hence, the sum of arguments  $p(u)$ and $p(w)$ cannot be equal to $\pi$ (-$\pi$), and therefore,  the equality  $p(u)p(w) = -c$ is impossible. Hence, $u$ and $w$ are both in the white zone. Let $u$ have the argument $\beta$. It follows that  $$|p(u)p(w)| = (1 + 2\cos(\beta))(1 + 2\cos(\pi - \beta)) = 1 - 4 \cos^2 \beta$$

For $\beta$ such that  $\cos \beta \ne 0$, this value is strictly less than one,  which is required. This condition can be violated only if either $u = \pm i$ or $w = \pm i$. 
Let  $u = i$. In this case  $p(u) = 1 + i - 1 = i$, and the product of the  numbers  $p(u)$ and  $p(w)$ can give a number with the argument  $\pi$ only if $p(w) = i$. Then $w = i$, and hence  $z^m = i$, $z^q = i$. 

Let us now prove that the assertion $m \not \equiv q \pmod{4}$ is impossible. 
Denote by  $\gamma$ the argument of $z$, chosen in the interval $[0, 2\pi)$. 
Since  $z^m = i$, we have $\gamma = \frac{\pi}{2m} + \frac{2 \pi m_1}{m}$, where $m_1 = 0, \ldots, m - 1$. Since  $z^q = i$, another formula is true as well:  $\gamma = \frac{\pi}{2q} + \frac{2 \pi q_1}{q}$, where $q_1 = 0, \ldots, q - 1$. Thus,  $\frac{1}{2m} + \frac{2 m_1}{m} = \frac{1}{2q} + \frac{2 q_1}{q}$, and after multiplication by  $2mq$ we obtain $q - m = 4 (q_1m - m_1q)$. This contradicts to our assumption.  

The case $u = -i$ is considered similarly. In this case  $w = -i$ and the number 
 $q - m$ must be again a multiple of $4$. 

{\hfill $\Box$}
\smallskip

{\tt Proof for series \ref{ser6}}.   
We begin with the following assertion: 

Consider a polynomial $P(z)  =  (1+z^m)q(z^r) +1$,  where  $q(z) = z^n + a_{n-1}z^{n-1} + \ldots + a_1z + 1$ is an integer polynomial. The polynomial $P$ is expanding if 
  $$\min \limits_{t \in \re} (\cos nt  + a_{n-1} \cos (n-1)t + \ldots + a_1 \cos t + 1) > - \frac{1}{2}.$$
 
Similarly to the previous series, we apply the maximum modulus principle and 
reduce the proof to the following statement: if $z = e^{i\alpha}$ and 
$(1 + z^m)q(z) = - c$, where $c$ is real and strictly positive, then $c < 1$. 
Note that the number  $1+z^m$ is located in the right half-plane (if $1 + z^m = 0$, then
 $c = 0$ is not positive). Since the product of $1+z^m$ and $q(z^r)$ has the 
 argument $\pi$, then $q(z^r)$ has a negative real part. Suppose $q(z^r)$ has the argument $\gamma$. Then $(1+z^m)$ has the argument $\pi - \gamma$, and $|1 + z^m| = |2\cos (\pi - \gamma)|$.  Therefore,  
\begin{multline*}
|(1+z^m)q(z^r)|  = 2|\cos (\pi - \alpha) q(z^r)|  = 2|\cos (\alpha) q(z^r)|   =  2 |\Re(q(z^r))| = \\
= -2 \Re(q(z^r)) = 2 (\cos nrs + a_{n-1} \cos (n-1)rs + \ldots + a_1 \cos rs + 1) < 2 \frac{1}{2} = 1 
\end{multline*}
Hence, $c< 1$, which concludes the proof. 

In our case  $q(z) = z^2 \pm z + 1$, $\min \limits_{t \in \R}(\cos 2t \pm \cos t +1)  = \min \limits_{x\in [-1, 1]} 2x^2 \pm x = - \frac{1}{8}$.

{\hfill $\Box$}
\smallskip

\begin{figure}[ht]
\centering
\includegraphics[width = 0.4\textwidth]{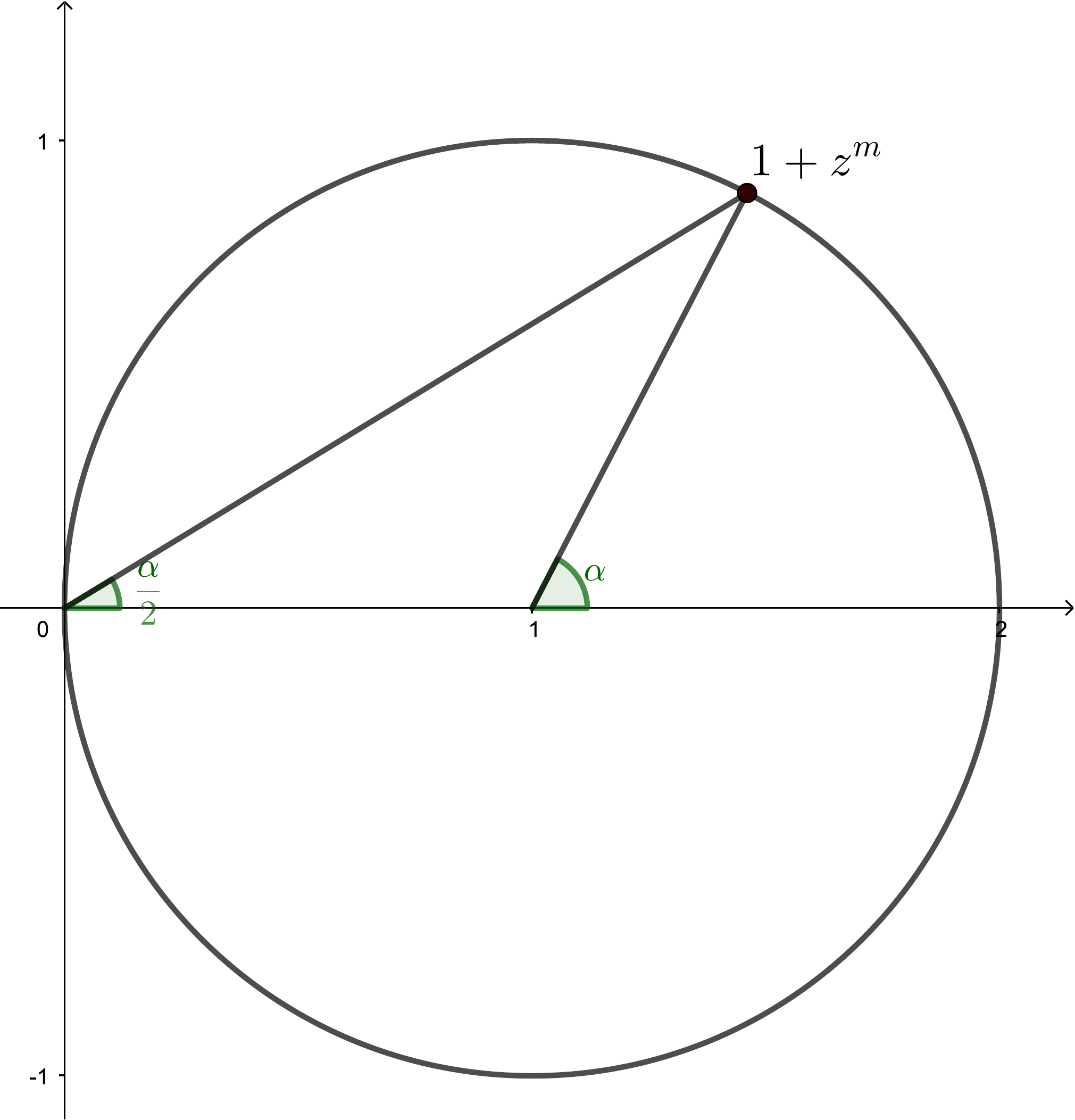}
\caption{The argument of the form $1 + z^m$}
\label{double}
\end{figure}

{\tt Proof of the series \ref{ser7}}.   
We apply the maximum modulus principle as in the previous two series. 
It suffices to prove that if $|z| = 1$ and  $P(z) - 1 = -c$, where $c$ is a real positive number, then $c < 1$. We prove that the equality $P(z) - 1 = -c$ does not hold 
for any  $c > 0$. Assume the converse. Let $s$ be the argument of the number  $z$; 
then the arguments of $1 + z^a$, $1 + z^b$, $1 + z^{(a + b)}$, $1 + z^{2 (a + b)}$, $\ldots$, $1 + z^{2^k (a + b)}$, up to addition of $\pi$, are equal to 
 $\frac{as}{2}$, $\frac{bs}{2}$, $\frac{(a + b)s}{2}$, $(a + b)s$, $\ldots$, $2^{k - 1} (a + b)s$ respectively, see Fig.~\ref{double}. 

Their sum, up to an addition of $2\pi$, is equal to $\pi$, hence for a suitable number 
$t \in \z$, one has  

$$\frac{as}{2} + \frac{bs}{2} + \frac{(a + b)s}{2} + (a + b)s + \ldots + 2^{k - 1}(a + b)s = \pi t.$$  

Hence, $\pi t = (a+b)s (1 + 1 + 2 + 4 + \ldots + 2^{k-1}) = (a+b)s 2^k$. 
Then the argument of the number   $z^{2^k (a + b)}$ in the last brackets is equal, up to 
an addition of $2\pi$, to  $\pi t$. This means that either $1 + z^{2^k (a + b)} = 0$ (and so  $P(z) - 1 = 0$, $c = 0$, which is a contradiction) or  $1 +  z^{2^k (a + b)} = 2$. 
Hence, now we can consider the inequality with a less number of brackets: 
$$(1 + z^a)(1 + z^b)(1 + z^{(a + b)}) (1 + z^{2(a + b)}) (1 + z^{4(a + b)}) (1 + z^{8 (a + b)}) \ldots (1 + z^{2^{k-1} (a + b)}) = -\frac{c}{2}.$$ 
Then we repeat this argument reducing the number of brackets. 
This way we come to the polynomial $(1 + z^a)(1 + z^b)$, 
which takes values only from the right half-plane 
(see the proof for series \ref{ser3}) and hence cannot take a real negative value.

{\hfill $\Box$}
\smallskip

Now we can give the promised proof of Theorem~\ref{t.ser4} from Section~10.  
 By $b_{+}(d)$ we denote the number of good partitions of the number $d$, 
 by $b_{-}(d)$ we denote the number of bad partitions (see Section~10 for the definitions). 
 By  $b(d)$ we denote the number of all partitions of the number $d$ 
 to the sum of natural terms. We have 
$b(d) = b_{-}(d)+b_{+}(d)$. 
We use the following well-known estimate:
\begin{lemma}\label{l.3-part}
For every $d \in \n$, we have    
$$ 
\frac{d(d-1) - r(r-1)}{12}\ -\ \frac{1}{2}\quad \le \quad b(d) \quad \le \quad  
\frac{d(d-1) - r(r-1)}{12}\ +\ \frac{1}{2}, 
$$
where $r$ is the remainder of~$d$ upon the division by $3$. 
\end{lemma}

Since the term $r(r-1)$ on the set of remainders of the division by~$3$
takes only the values~$0$ and $2$, we obtain 
\begin{cor}\label{c.3-part}
For every $d \in \n$, we have   
$\ b(d) \ = \ \frac{d(d-1)}{12} \ + \ \omega_d\, ,$ 
where $\ \omega_d  \in \bigl[-\frac23, \frac12 \bigr]\, .$ 
\end{cor}

{\tt Proof of Theorem~\ref{t.ser4}}.   
The sum of components of a bad vector is divisible by  $3$, 
hence, for $d \not \equiv 0 \pmod 3$, there are no bad partitions. 
Therefore,~$b_+(d) \, = \, b(d)$, and we apply Corollary~\ref{c.3-part} to this value. 
If, on the other hand, 
$d  \equiv 0 \pmod 3$, then $\, d \, = \, 3^{\ell}\, d', \ d' \not \equiv 0 \pmod 3$. If  $d=n_1+n_2+n_3$ is a bad partition, then there exists  $j\le \ell-1$ for which 
 $(n_1, n_2, n_3) =  3^{j}(3x+s, 3y+s, 3z+s)$, where  $s\in \{1,2\}$
and $x,y,z \in \z$. 
Consequently, $x+ y+ z\,  = \, 3^{-j-1}d \ - \, s$. Thus, every bad partition of the number  $d$ has the form~$d \, = \, 3^{j}(3x+s) + 3^{j}(3y+s)+ 3^{j}(3z+s)$, 
where $j= 0, \ldots , \ell-1, \, s \in \{1, 2\}$ and  $x+y+z$ is an arbitrary partition of a number  $3^{-j-1}d \ - \, s$. Therefore, if  $b\equiv 0\pmod 3$,  then  
\begin{equation}\label{eq.sum-quadr}
b_{+}(d) \ = \ b(d) \ - \ \sum_{j=0}^{\ell-1} \, \left(b\left(3^{-j-1}d \ - \, 1 \right)\ + \ 
b\left(3^{-j-1}d \ - \, 2 \right)\right)
\end{equation}
Denote $a_j = 3^{-j-1}d$. All the numbers $a_j$
except for  $a_{\ell-1}$ are multiples of $3$.
Hence, by Lemma~\ref{l.3-part}, the $j$th term in the sum~(\ref{eq.sum-quadr}) 
for $j = 0, \ldots , \ell-1$,  is equal to 
$$
\frac{(a_j-1)(a_j-2) - 2\cdot 1}{12} \ + \ \frac{(a_j-2)(a_j-3)}{12} \ +\ 
v_j \quad = \quad   \frac{(a_j - 1)(a_j-3)}{6}  \ +\ 
v_j, 
$$
where $v_j \in [-1,1]$. In the latter term 
$a_{\ell -1}$ is not divisible by 3, hence, the numerator of the first fraction 
can increase by $2$. Thus, 
$$
b_+(d) \ = \ \frac{d(d-1)}{12} \ - \ v \, - \ \sum_{j=0}^{\ell-1} \, \left(\frac{(a_j - 1)(a_j-3)}{6}  \ +\ 
v_j\right) \, , 
$$
where $v \in \bigl[-\frac13,\frac56\bigr], \, $ and $\ v_j \in \bigl[-1, 1 \bigr]$. 
To estimate this value from below we replace the sum of the terms $\frac{a_j^2}{6}$ by an infinite sum and the sum $\frac{2a_j}{3}$ by the first term. Taking the sum of the geometric progression and having estimated  $\sum \limits_{j=0}^{l-1} (0.5 + v_j) \le 1.5 l \le 1.5 d$, we obtain the required lower bound. To get an upper bound we on the contrary replace the sum of  $\frac{a_j^2}{6}$ by the first term and the sum $\frac{2a_j}{3}$ by the infinite geometric progression. Thus, we obtain the upper bound. 

{\hfill $\Box$}
\smallskip

\bigskip

\noindent \textbf{Acknowledgements.} The authors are grateful for useful discussions to  S.Konyagin, A.Dubickas, I.Shkredov, N.Moshchevitin, and S.Ivanov. 

The work of the second author is supported by the  Foundation for Advancement of Theoretical Physics and Mathematics ``BASIS''.   
\bigskip

\bigskip

\begin{center}
\textbf{Appendix}
\end{center}
\bigskip

{\tt Proof of Lemma~\ref{l.dyn}}. Assume $\cT(M_1, \bar \ba_1) = \cT(M_2, \bar \ba_2)$. 
First we show that  $\rho(M_1^{-1}) = \rho(M_2^{-1})$, where $\rho$ 
is the spectral radius of the matrix.  Since $M_1$ does not have multiple eigenvalues  
and since  $\ba_1$ does not belong to an invariant subspace of~$M_1$, we obtain 
$\|M_1^{-k} \ba_1\|\, \asymp \,  [\rho(M_1^{-1})]^k$, where 
$\asymp$ denotes the asymptotic equivalence. Then the number of segments in the set  $\cT(M_1, \bar \ba_1)$ whose length exceeds a given number $\varepsilon$ is equal to   
 $\frac{\ln \varepsilon }{\ln \rho(M_1^{-1})} \, + \, o\bigl(1\bigr)$
as $\varepsilon \to 0$. Since  $\cT(M_1, \bar \ba_1) = \cT(M_2, \bar \ba_2)$, 
we obtain $\rho(M_1^{-1})= \rho(M_2^{-1})$. 

It may be  assumed that~$M_1^{-1}$ has a unique largest by modulus eigenvalue 
(or two complex conjugate), otherwise the problem is reduced to several problems of smaller dimensions.  We call the largest by modulus eigenvalue 
{\em leading}, and its eigenvector also call the {\em leading} eigenvector. Let
 $V_i$ be the subspace in~$\re^d$ spanned by the leading eigenvectors of the matrix~$M_i^{-1}$,  
  $i = 1,2$. Two cases are possible. 
  
  \smallskip 
  
 \textbf{1)}  ${\dim V_1} = 1$. In this case the leading eigenvalue 
 $\lambda$ is real and  $V_1$ is spanned by the leading eigenvector $\bv_1$. 
 Then the direction of the vector $M_1^{-k}\bar \ba_1$ converges to the direction of  
 $V_1$ as $k \to \infty$. Hence, the ratio of lengths of the successive segments  
 $M_1^{-k-1}\bar \ba_1$  and $M_1^{-k}\bar \ba_1$ tends to 
 $\rho(M_1^{-1})$ 
as $k \to \infty$. Since $\rho(M_1^{-1}) < 1$, it follows that for sufficiently large $k$, 
we know the order of the elements in the sequence  $\{M_1^{-k}\bar \ba_1\}_{k\ge N}$. 
 This means that for small~$\varepsilon > 0$, we know this order of elements 
 of lengths smaller than~$\varepsilon$ in the set~$\cT(M_1 , \bar \ba_1)$. 
 We take  $d+2$ segments in a row and denote by  $\ell_1, \ldots \ell_{d+2}$
the straight lines on which they lie. 
The operator $M_1^{-1}$ maps the line $\ell_i$ to $\ell_{i+1}, \, i = 1, \ldots , d+1$. 
By considering each line as a point in the projective space  
  $\mathbb{P}^{d-1}$ we see that the operator $M_1^{-1}$ is uniquely defined, up to a sign, 
  as an operator in~$\mathbb{P}^{d-1}$ that maps given $d+1$ points to given points. 
  Now consider the operator~$M_2^{-1}$. It must have the same leading eigenvector~$\bv_2 = \bv_1$ as a unique limit direction of vectors from the set~$\cT(M_2, \bar \ba_2)$
  (hence,~${\dim V_2} = 1$). It must have the same order of segments of lengths smaller than~$\varepsilon$ and therefore, it has the same lines
   $\ell_1, \ldots , \ell_{d+2}$ in the same order. This implies that  $M_1$ and $M_2$
   define the same projective operator in the space $\mathbb{P}^{d-1}$, 
   hence $M_1 = \pm M_2$. 
\smallskip 

\textbf{2)}  ${\dim V_1} = 2$. In this case the leading eigenvalue of  $M_1^{-1}$
has the form $r \, e^{\, 2\pi i \alpha_1}$, where $r = \rho(M_1^{-1})$. We need to consider two cases.

\textbf{a)} $\, \alpha_1 = \frac{p}{q}$ is rational ($\frac{p}{q}$ is an irreducible 
fraction). In a suitable basis in the two-dimensional plane~$V_1$ the operator  $M_1^{-1}|_{V_1}$ is a rotation by the angle $\frac{2\pi p}{q} $ with a multiplication by $r$. 
The segments $M_1^{-k}\bar \ba_1$ have exactly  $q$ limit directions
corresponding to the angles~$\frac{2\pi n}{q}, \, n= 0, \ldots , q-1$. There are at least 
three these directions, because the case~$q=2$ corresponds to the case 1) of real $\lambda$. 
Therefore, there exists a unique linear transform that maps the limit directions of 
vectors of the set $\cT(M_1, \bar \ba_1)$ to the direction of the angles~$\frac{2\pi p}{q}, \, p= 0, \ldots , q-1$. Consequently, 
the basis on the plane~$V_1$ in which the operator $M_1^{-1}$ is a rotation with a contraction 
is uniquely defined by the set~$\cT(M_1, \bar \ba_1)$. In this basis 
the ratio of lengths of consecutive segments in the set $\{M_1^{-k}\bar \ba_1\}_{k \ge N} $  
tends to~$r$, which gives us the order of elements of this sequence. 
The following reasoning is the same as in case~1): having obtained the order of the sequence 
we restore the operator~$M_1^{-1}$ and thus prove that 
$M_1 = \pm M_2$.   
\smallskip 

\textbf{b)} $\, \alpha_1$ is irrational. In this case, in a suitable basis on the 
two-dimensional plane~$V_1$, the operator $M_1^{-1}|_{V_1}$ is a rotation by an 
irrational angle $2\pi \alpha_1$ with a multiplication by $r$. 
To identify that basis we apply the 
Weyl theorem on uniform distribution~\cite{Weil}:
a trajectory of an irrational rotation of a unit  vector converges to a uniform distribution on a circle. This means that for every arc of a unit circle, the number 
of points came to that arc to the total number of points of the sequence
from the first one to the $n$th one tends to the ratio of the arc length to the 
length of the circle. The number of segments in the set 
$\cT(M_1, \bar \ba_1)$ of length larger than $\varepsilon$ is equal to 
 $\frac{\ln \varepsilon }{\ln r} \, + \, O\bigl(1\bigr)$
as $\varepsilon \to 0$. The number of those segments in a given angle of the value $\beta$ is equal to  $\frac{\beta\ln \varepsilon }{\pi \ln r} \, + \, o\bigl(1\bigr)\, , \ \varepsilon \to 0$.  
If there exists another suitable basis in~$V_1$, then it has to define the same value 
of the angle~$\beta$. Thus, in the new basis the values of all angles are the same as in the old one. Hence, this is the same basis, up to an orthogonal transform.  
We see that  {\em for an arbitrary set $\cT(M_1, \bar \ba_1)$, 
 there exists a unique, up to an orthogonal transform, basis in the two-dimensional 
 plane~$V_1$ in which the operator  $M_1^{-1}$ is a rotation with a multiplication by~$r$.} 
 In this basis, the order of the segments in this sequence is uniquely defined
 and hence, the operator~$M_1^{-1}$ is uniquely defined as well, up to a sign. 
If we consider now the operator~$M_2^{-1}$ with the same set~$\cT(M_2, \bar \ba_2)= \cT(M_1, \bar \ba_1)$, then it will have the same two-dimensional plane~$V_2 = V_1$
as a limit plane of directions of segments from~$\cT(M_2, \bar \ba_2)$
with the same basis in which~$M_2^{-1}$ defines a rotation with a multiplication. 
Therefore, the sequence $M_2^{-k}\bar \ba_2$ has the same order of elements, provided 
$k$ is large enough.  In the same way as in item~1), 
having learned the order of elements we recover the operator~$M_2^{-1}$ 
and thus prove that 
$M_1 = \pm M_2$.

{\hfill $\Box$}
\smallskip

{\tt Proof of Lemma~\ref{l.3-part}}. 
Without loss of generality we assume that  $n_1 \le n_2 \le n_3$. 
Let  $d = 3k + r, \, r \in \{0, 1, 2\}$. 
For $k=0$, the assertion is obvious, hence, we suppose $k \ge 1$. 
Then  $n_1 \le k$, and $n_2 \in \bigl[n_1 , \bigl[\frac{d-n_1}{2}\bigr]\,  \bigr]$, 
where $[x]$ is an integer part of $x$. Consequently, 
$$
b(d) \ = \ \sum_{n_1 = 1}^{k}\, \left(\left[\frac{d-n_1}{2}\right]\, - \, n_1 \, +\, 1\right)\, . 
$$
The total number of odd numbers in the sequence   $d-n_1\, , \, n_1 = 1, \ldots , k$, 
is at least  $\frac{k-1}{2}$ and at most  $\frac{k+1}{2}$. 
Hence,  
$$
- \, \frac{k+1}{2} \, +\, \sum_{n_1 = 1}^{k}\, \left(\frac{d-n_1}{2}\, - \, n_1 \, +\, 1\right)\quad  \le \quad b(d) \quad  \le \quad  - \, \frac{k-1}{2}\, + \, \sum_{n_1 = 1}^{k}\, \left(\frac{d-n_1}{2}\, - \, n_1 \, +\, 1\right)\, .
$$
Now by a direct calculation we complete the proof.  

{\hfill $\Box$}


\smallskip


\begin{thebibliography}{NN}
\bibitem[ABPT]{ABPT}
S.\,Akiyama, H.\,Brunotte, A.\,Peth\H o, J.\,Thuswaldner, 
\newblock {\em Generalized radix representations and dynamical systems III},  
\newblock Osaka J. Math. 45 (2008), no. 2, 347 -- 374.
\smallskip


\bibitem[AG]{AG}
S.\,Akiyama, N.\,Gjini,
\newblock {\em On the connectedness of self-affine attractors}, 
\newblock Archiv der Mathematik,  82 (2004), no. 2,  153 -- 163. 
\smallskip


\bibitem[AL]{AL}
S.\,Akiyama, B.\,Loridant,  
\newblock {\em Boundary parametrization of planar self-affine tiles with collinear digit set},  
\newblock Science China Mathematics 53 (2010), no. 9, 2173 -- 2194. 
\smallskip


\bibitem[ALT]{ALT}
S.\,Akiyama, B.\,Loridant, J.M.\,Thuswaldner,
\newblock {\em Topology of planar self-affine tiles with collinear digit set},  
\newblock (2018) arXiv:1801.02957   
\smallskip

\bibitem[AP]{AP}
S.\,Akiyama, S.A.A.\,Peth\H o, 
\newblock {\em On the distribution of polynomials with bounded roots II. Polynomials with integer coefficients},  
\newblock Uniform Distribution Theory 9 (2014), no. 1, 5 -- 19.
\smallskip



\bibitem[AT]{AT}
S.\,Akiyama, J.\,Thuswaldner,
\newblock {\em A survey on topological properties of tiles related to number systems},  
\newblock Geom Dedicata, 109 (2004), no. 1, 89 -- 105. 
\smallskip


\bibitem[B10]{B10}
C.\,Bandt,
\newblock {\em Combinatorial topology of three-dimensional self-affine tiles},  
\newblock (2010) arXiv:1002.0710 
\smallskip

\bibitem[B91]{B91}
C.\,Bandt,
\newblock {\em  Self-similar sets. V. Integer matrices and fractal tilings of $\re^n$},  
\newblock Proc. Amer. Math. Soc. 112 (1991), no. 2, 549 -- 562.
\smallskip


\bibitem[BG]{BG}
C.\,Bandt, G.\,Gelbrich, 
\newblock {\em  Classiffication of self-affine lattice tilings},  
\newblock J. London Math. Soc. 50 (1994), no. 3, 581 -- 593.
\smallskip


\bibitem[BL]{BL}
J.J.\,Benedetto, M.\,Leon, 
\newblock {\em The construction of multiple dyadic minimally supported frequency wavelets on $\re^d$},  
\newblock Contemp. Math. 247 (1999), 43 -- 74.
\smallskip

\bibitem[BS]{BS}
J.J.\,Benedetto, S.\,Sumetkijakan, 
\newblock {\em Tight frames and geometric properties of wavelet sets},  
\newblock Adv. Comput. Math. 24 (2006), no. 1--4, 35 -- 56.
\smallskip


\bibitem[BW]{BW}
C.\,Bandt, Y.\,Wang, 
\newblock {\em  Disk-Like Self-Affine Tiles in $\re^2$},  
\newblock Discrete and Computational Geometry 26 (2001), no. 4, 591 -- 601.
\smallskip

\bibitem[Bow]{Bow}
M.\,Bownik, 
\newblock {\em  Anisotropic Hardy spaces and wavelets}, 
\newblock  Mem. Amer. Math. Soc. 164 (2003), no. 781.
\smallskip



\bibitem[CDM]{CDM}
A.S.\,Cavaretta, W.\,Dahmen~W, and C.A.\,Micchelli, 
\newblock {\em Stationary subdivision}, 
\newblock  Mem. Amer. Math. Soc. 93 (1991), no. 453. 
\smallskip 

\bibitem[CGRS]{CGRS}
M.\,Cotronei, D.\,Ghisi, M.\,Rossini, T.\,Sauer, 
\newblock {\em An anisotropic directional subdivision and multiresolution scheme}, 
\newblock  Adv. Comput. Math. 41 (2015), no. 3, 709 -- 726. 
\smallskip 

\bibitem[CHM]{CHM}
C.A.\,Cabrelli, C.\,Heil, U.M.\,Molter,
\newblock {\em Self-similarity and multiwavelets in higher dimensions},
\newblock Memoirs Amer. Math. Soc. 170 (2004), no. 807.
\smallskip

\bibitem[CM]{CM}
M.\,Charina, T. Mejstrik, 
\newblock {\em Multiple multivariate subdivision schemes: Matrix and operator approaches},  \newblock Comp. Appl. Math., 349 (2019), 279 -- 291.
\smallskip 

\bibitem[CP]{CP}
M.\,Charina, V.Yu.\,Protasov,  
\newblock {\em Regularity of anisotropic refinable functions},  
\newblock Appl. Comput. Harmon. Anal., 47 (2019), no. 3, 795 -- 821.
\smallskip

\smallskip
\bibitem[CGV]{CGV}
A.\,Cohen, K.\,Gr\"ochenig, and L.F.\,Villemoes, 
\newblock {\em Regularity of multivariate refinable functions}, 
\newblock Constr. Approx. 15 (1999), no. 2, 241 -- 255.
\smallskip

\bibitem[CT]{CT}
G.R.\,Conner, J.M.\,Thuswaldner, 
\newblock {\em Self-affine manifolds},  
\newblock Advances in Mathematics 100 (2016), no. 289, 725 -- 783.
\smallskip

\bibitem[Desp]{Desp}
N.\,Desprez, 
\newblock {\em Chaoscope software},  
\newblock \href{http://www.chaoscope.org/}{http://www.chaoscope.org/}
\smallskip

\bibitem[DK]{DK}
A.\,Dubickas,  S.V.\,Konyagin, 
\newblock {\em On the number of polynomials of bounded measure}, 
\newblock Acta Arith. 86 (1998), no. 4, 325 -- 342.
\smallskip 

\bibitem[DL]{DL}
Q.\,Deng, K.-S.\,Lau,  
\newblock {\em Connectedness of a class of planar self-affine tiles},  
\newblock Journal of Mathematical Analysis and Applications 380 (2011), no. 2, 493 -- 500.
\smallskip

\bibitem[DLS]{DLS}
X.\,Dai, D.R.\,Larson, D.M.\,Speegle,
\newblock {\em Wavelet sets in $\re^n$},  
\newblock J. Fourier Anal. Appl. 3 (1997), no. 4, 451 -- 456.
\smallskip



\bibitem[Dub]{Dub}
A.\,Dubickas, 
\newblock {\em Counting integer reducible polynomials with bounded measure}, 
 Appl. Anal. Discrete Math. 10 (2016), no. 2, 308 -- 324.
\smallskip 


\bibitem[FG]{FG}
X.\,Fu, J.-P.\,Gabardo, 
\newblock {\em Self-affine scaling sets in $\re^2$},  
\newblock Memoirs of the American Mathematical Society 233 (2015), 1 -- 97.
\smallskip


\bibitem[G81]{G81}
W.J.\,Gilbert, 
\newblock {\em Radix representations of quadratic fields},  
\newblock J. Math. Anal. Appl. 83 (1981), no. 1, 264 -- 274.
\smallskip



\bibitem[Gar]{Gar}
A.\,Garsia, 
\newblock {\em Arithmetic properties of Bernoulli convolutions},
\newblock Trans. Amer. Math. Soc. 102 (1962), no. 3, 409 -- 432.
\smallskip



\bibitem[Gel]{Gel}
G.\,Gelbrich, 
\newblock {\em Self-affine Lattice Reptiles with Two Pieces in $\re^n$},  
\newblock Math. Nachr., 178 (1996), no. 1, 129 -- 134.
\smallskip

\bibitem[GJ]{GJ}
R.\,Gundy, A.\,Jonsson, 
\newblock {\em Scaling functions on $\re^2$ for dilations of determinant $\pm 2$},  
\newblock Appl.  Comput. Harmon. Anal. 29 (2010), no. 1, 49 -- 62.
\smallskip




\bibitem[GM]{GM}
C.\,Gr\"ochenig, W.R.\,Madych, 
\newblock {\em  Multiresolution analysis, Haar bases, and self-similar tilings of $\re^n$},  
\newblock IEEE Trans. Inform. Theory 38 (1992), no. 2, 556 -- 568.
\smallskip

\bibitem[GH]{GH}
K.\,Gr\"ochenig, A.\,Haas,
\newblock {\em Self-similar lattice tilings},
\newblock J. Fourier Anal. Appl., 1 (1994), no. 2, 131 -- 170.
\smallskip

\bibitem[HL]{HL}
X.G.\,He, K.S.\,Lau,  
\newblock {\em Characterization of tile digit sets with prime determinants},  
\newblock Appl. Comput. Harmonic Anal., 16 (2004), no. 3, 159 -- 173.
\smallskip

\bibitem[HSV]{HSV}
D.\,Hacon, N.C.\,Saldanha, J.J.P.\,Veerman, 
\newblock {\em Remarks on self-affine tilings},  
\newblock Exp. Math. 3 (1994), no. 4, 317 -- 327.
\smallskip

\bibitem[KL00]{KL00}
I.\,Kirat, K.-S.\,Lau, 
\newblock {\em On the connectedness of self-affine tiles},  
\newblock J. Lond. Math. Soc., 62 (2000), no. 1, 291 -- 304.
\smallskip


\bibitem[KL02]{KL02}
I.\,Kirat, K.-S.\,Lau, 
\newblock {\em Classification of integral expanding matrices and self-affine tiles},  
\newblock Discrete Comput. Geom. 28 (2002), no. 1, 49 -- 73.
\smallskip


\bibitem[KLR]{KLR}
I.\,Kirat, K.-S.\,Lau, H.\,Rao, 
\newblock {\em Expanding polynomials and connectedness of self-affine tiles},  
\newblock Discrete Comput. Geom. 31 (2004), no. 2, 275 -- 286.
\smallskip


\bibitem[KM]{KM}
A.\,Kravchenko, D.\,Mekhontsev, 
\newblock {\em IFS Builder 3d software},  
\newblock \href{http://fractals.nsu.ru/builder3d_en.htm}{http://fractals.nsu.ru/builder3d_en.htm}
\smallskip

\bibitem[KPS]{KPS}
A.~Krivoshein, V.\,Yu.~Protasov, and M.\,A.~Skopina,
\newblock {\em Multivariate wavelets frames}, 
\newblock Springer, 2016
\smallskip 


\bibitem[KT]{KT}
P.\,Kirschenhofer, J.M.\,Thuswaldner, 
\newblock {\em Shift radix systems -- a survey},  
\newblock In: Numeration and Substitution 2012, Research Institute for Mathematical Sciences (RIMS), Kyoto (2014), 1 -- 59.
\smallskip

\bibitem[KT17]{KT17}
P.\,Kirschenhofer, A.\,Thuswaldner, 
\newblock {\em Distribution results on polynomials with bounded roots},  
\newblock Monatshefte f\"ur Mathematik 185 (2018), no. 4, 689 -- 715.
\smallskip


\bibitem[LW95]{LW95}
J. C.\,Lagarias, Y.\,Wang,  
\newblock {\em Haar type orthonormal wavelet bases in $\re^2$},  
\newblock J. Fourier Anal. Appl. 2 (1995), no. 1, 1 -- 14.
\smallskip


\bibitem[LW96]{LW96}
J.\,Lagarias, Y.\,Wang, 
\newblock {\em Haar bases for $L_2(\re^n)$ and algebraic number theory},   
\newblock J. Number Theory 57 (1996), no. 1, 181 -- 197. 
\smallskip

\bibitem[LW97]{LW97}
J.\,Lagarias, Y.\,Wang, 
\newblock {\em Integral self-affine tiles in $\re^n$. II. Lattice tilings},  
\newblock J. Fourier Anal. Appl. 3 (1997), no. 1, 83 -- 102.
\smallskip


\bibitem[LW99]{LW99}
J.\,Lagarias, Y.\,Wang, 
\newblock {\em Corrigendum/addendum: Haar bases for $L_2(\re^n)$ and algebraic number theory, J. Number Theory 57 (1996), no. 1, 181 -- 197},  
\newblock  J. Number Theory 76 (1999), no. 2, 330 -- 336. 
\smallskip



\bibitem[Mej]{Mej}
T.\,Mejstrik, 
\newblock {\em t-toolboxes for Matlab}, 
\newblock  \href{tommsch.com/science.php}{tommsch.com/science.php}, 2018
\smallskip 


\bibitem[Mekh]{Mekh}
D.\,Mekhontsev, 
\newblock {\em IFStile software},  
\newblock \href{http://ifstile.com/view/Main_Page}{http://ifstile.com/view/Main_Page}
\smallskip


\bibitem[Mer15]{Mer15}
K.D.\,Merrill,  
\newblock {\em Simple Wavelet Sets in $\re^n$},  
\newblock J. Geom. Anal. 25 (2015), no. 2, 1295 -- 1305.
\smallskip

\bibitem[Mer18]{Mer18}
K.D.\,Merrill,  
\newblock {\em Generalized Multiresolution Analyses. Applied and Numerical Harmonic Analysis},  
\newblock Springer, Berlin (2018).
\smallskip

\bibitem[Mor]{Mor}
F.\,Morgan, 
\newblock {\em Geometric Measure Theory}, 
\newblock   Academic Press (1988).
\smallskip

\bibitem[NPS]{NPS}
I.\,Novikov, V.Yu.\,Protasov, M.A.\,Skopina,
\newblock {\em Wavelets theory}, 
\newblock AMS, Translations Mathematical Monographs, 239 (2011).
\newblock  


\bibitem[NSVW]{NSVW}
S.-M.\,Ngai, V.F.\,Sirvent, J.J.P.\,Veerman, Y.\,Wang, 
\newblock {\em On 2-reptiles in the plane},  
\newblock Geom. Dedicata 82 (2000), no. 1, 325 -- 344.
\smallskip

\bibitem[P97]{P97}
V.\,Yu.~Protasov,
\newblock {\em The generalized spectral radius. A geometric approach},
\newblock Izvestiya Math., 61 (1997), 995 -- 1030.
\smallskip


\bibitem[P20]{P20}
V.Yu.\,Protasov, 
\newblock {\em Surface dimension,  tiles, and synchronising automata},  
\newblock SIAM J. Math. Anal., to appear (2020). 
\smallskip

\bibitem[Rad]{Rad}
H.\,Rademacher, 
\newblock {\"Uber partielle und totale Differenzierbarkeit I.} 
\newblock Math. Ann. 89 (1919), 340 -- 359. 
\smallskip

\bibitem[ST]{ST}
W.\,Steiner, J.M.\,Thuswaldner, 
\newblock {\em Rational self-affine tiles},  
\newblock Trans. Amer. Math. Soc. 367 (2015), no. 11, 7863 -- 7894.
\smallskip


\bibitem[TZ]{TZ}
J.M.\,Thuswaldner, S.-Q.\,Zhang, 
\newblock {\em On self-affine tiles whose boundary is a sphere},  
\newblock Trans. Amer. Math. Soc. 373 (2020), no. 1, 491 -- 527.
\smallskip


\bibitem[U19]{U19}
M.\,Uray,
\newblock {\em On the expansivity gap of integer polynomials},  
\newblock (2019) arXiv:1905.06976 
\smallskip

\bibitem[Weil]{Weil}
H.\,Weyl, 
\newblock{ \"Uber die Gleichverteilung von Zahlen mod. Eins}, 
\newblock Mathematische Annalen, 77 (1916), 313 -- 352.
\smallskip 


\bibitem[Woj]{Woj}
P.\,Wojtaszczyk,  
\newblock {\em A Mathematical Introduction to Wavelets},  
\newblock London Math. Soc. Stud. Texts, vol. 37, Cambridge Univ. Press, Cambridge, New York, Melbourne, Madrid (1997). 
\smallskip

\bibitem[Zai]{Zai}
T.\,Zaitseva, 
\newblock {\em Haar wavelets and subdivision algorithms on the plane},  
\newblock Advances in Systems Science and Applications 17 (2017), no. 3,  49 -- 57.
\smallskip


\bibitem[Zai2]{Zai2}
T.\,Zaitseva, 
\newblock {\em Simple tiles and attractors}, 
\newblock Sb. Math. 211 (2020), no 9. 
\smallskip 


\bibitem[Zakh]{Zakh}
V.G.\,Zakharov, 
\newblock {\em Rotation properties of 2D isotropic dilation matrices},  
\newblock Int. J. Wavelets Multiresolut. Inf. Process. 16 (2018), no. 01.
\smallskip


\end{thebibliography}
\end{document}